\documentclass[11pt,letterpaper]{article}  

\usepackage{float}
\usepackage{mathrsfs}
\usepackage{cite}
\usepackage[T1]{fontenc}
\usepackage[latin9]{inputenc}
\usepackage{amsmath}
\usepackage{amssymb}
\usepackage{dsfont}
\usepackage{graphicx}
\usepackage{enumerate}
\usepackage{color}
\usepackage{float}
\usepackage[]{algorithm2e}
\SetKwInOut{KwInitialize}{initialize:} 
\SetKwBlock{KwGenGraph}{co-design communication graph}{end} 
\SetKwBlock{KwGetK}{refine optimal controller}{end} 
\SetKwInOut{Input}{input}\SetKwInOut{Output}{output}\SetKwInOut{Return}{return}

\usepackage[margin=1in]{geometry}    
\makeatletter

\pdfpageheight\paperheight
\pdfpagewidth\paperwidth

\floatstyle{ruled}
\newfloat{algorithm}{tbp}{loa}
\providecommand{\algorithmname}{Algorithm}
\floatname{algorithm}{\protect\algorithmname}

\@ifundefined{showcaptionsetup}{}{%
 \PassOptionsToPackage{caption=false}{subfig}}
\usepackage{subfig}
\makeatother

\def\QED{~\rule[-1pt]{5pt}{5pt}\par\medskip}

\newtheorem{rem}{Remark}
\newtheorem{lem}{Lemma}

\newtheorem{prop}{Proposition}

\newtheorem{thm}{Theorem}

\newtheorem{coro}{Corollary}
\newenvironment{IEEEproof}{{\bf Proof: }}{ \hfill \QED}
\newtheorem{example}{Example}
\newcommand{\N}{\mathscr{N}}

\newcommand{\s}{\mathcal{S}}

\newcommand{\Rp}{\mathcal{R}_p}
\newcommand{\Y}{\mathcal{F}}
\newcommand{\Htwo}{\mathcal{H}_2}
\newcommand{\RHtwo}{\mathcal R\Htwo}
\newcommand{\R}{\mathbb{R}}

\newcommand{\Z}{\mathbb{Z}}

\newcommand{\Hinf}{\mathcal{H}_\infty}
\newcommand{\RHinf}{\mathcal{R}\Hinf}
\newcommand{\A}{\mathscr{A}}

\newcommand{\bsupp}[1]{\text{bsupp}\left(#1\right)}
\newcommand{\supp}[1]{\text{bsupp}\left(#1\right)}

\newcommand{\approxcitations}{\cite{Bon91,Pis86}}

\newcommand{\minimize}[1]{\underset{#1}{\text{minimize}}}
\newcommand{\Htwonorm}[1]{\left\| #1 \right\|_{\Htwo}}
\newcommand{\commnorm}[1]{\left\|#1\right\|_{\text{comm}}}

\newcommand{\vct}[1]{\mathbf{#1}}
\newcommand{\Zp}{\mathbb Z_+}
\newcommand{\Zpp}{\mathbb Z_{++}}
\newcommand{\kcommnorm}[2]{\left\|#2\right\|_{#1-\text{comm}}}

\newcommand{\edit}[1]{\textcolor{black}{#1}}

\begin{document}

\title{Communication Delay Co-Design in $\Htwo$ Distributed Control Using Atomic Norm Minimization}

\author{Nikolai Matni
\thanks{N. Matni is with the Department of Control and Dynamical Systems, California Institute of Technology. 1200 E California Blvd., Pasadena, CA, 91125. (626) 395-6247.
 \tt{\small nmatni@caltech.edu}.}
\thanks{This research was in part supported by the NSF, AFOSR, ARPA-E, and the Institute for Collaborative Biotechnologies through grant W911NF-09-0001 from the U.S. Army Research Office. The content does not necessarily reflect the position or the policy of the Government, and no official endorsement should be inferred.  
A preliminary version of this work \cite{M_CDC13_codesign} has appeared at the 52nd Annual Conference on Decision and Control in December, 2013. }}
\maketitle
\begin{abstract}

When designing distributed controllers for large-scale systems, the actuation, sensing and communication architectures of the controller can no longer be taken as given.  \edit{In particular, controllers implemented using dense architectures typically outperform controllers implemented using simpler ones -- however, it is also desirable to minimize the cost of building the architecture used to implement a controller.}  The recently introduced Regularization for Design (RFD) framework poses the controller architecture/control law co-design problem as one of jointly optimizing the competing metrics of controller architecture cost and closed loop performance, and shows that this task can be accomplished by augmenting the variational solution to an optimal control problem with a suitable atomic norm penalty.  Although explicit constructions for atomic norms useful for the design of actuation, sensing and joint actuation/sensing architectures are introduced, no such construction is given for atomic norms used to design communication architectures.  This paper describes an atomic norm that can be used to design communication architectures for which the resulting distributed optimal controller is specified by the solution to a convex program.  Using this atomic norm we then show that in the context of $\Htwo$ distributed optimal control, the communication architecture/control law co-design task can be performed through the use of finite dimensional second order cone programming.

\end{abstract}

\section{Introduction}
\label{sec:introduction}
Large-scale systems represent an important class of application areas for the control engineer -- prominent examples include the smart-grid, software defined networking (SDN) and automated highways.  For such large-scale systems, designing the controller \emph{architecture} -- placing sensors and actuators as well as the communication links between them -- is now also an important part of the controller synthesis process. Indeed controllers with denser actuation, sensing and communication architectures will typically outperform those with simpler architectures -- \edit{however it is also desirable to minimize the cost of constructing a controller architecture.}

In \cite{MC_RFD_TAC}, the author of this paper and V. Chandrasekaran address the problem of jointly optimizing the architectural complexity of a distributed optimal controller and the closed loop performance that it achieves by introducing the Regularization for Design (RFD) framework.  In RFD, controllers with complicated architectures are viewed as being composed of atomic controllers with simpler architectures -- this family of simple controllers is then used to construct various \emph{atomic norms} \cite{CRPW12,Bon91,Pis86} that penalize the use of specific architectural resources, such as actuators, sensors or additional communication links.  These atomic norms are then added as a penalty function to the variational solution to an optimal control problem \edit{(formulated in the model matching framework)}, allowing the controller designer to explore the tradeoff between architectural complexity and closed loop performance by varying the weight on the atomic norm penalty in the resulting convex optimization problem.  

In \cite{MC_RFD_TAC} we give explicit constructions of atomic norms useful for the design of actuation, sensing and joint actuation/sensing architectures, but do not address how to construct an atomic norm for communication architecture design.  Indeed constructing a suitable atomic norm for communication architecture design has substantial technical challenges that do not arise in actuation and sensing architecture design: we address these challenges in this paper.
We model a distributed controller as a collection of sub-controllers, each equipped with a set of actuators and sensors, that exchange their respective measurements with each other subject to \emph{communication delays} imposed by an underlying communication graph.  Keeping with the philosophy adopted in RFD \cite{MC_RFD_TAC}, we view dense communication architectures, i.e., ones with a large number of communication links between sub-controllers, as being composed of multiple simple \emph{atomic} communication architectures, i.e., ones with a small number of communication links between sub-controllers. 
Thus the problem of controller communication architecture/control law co-design can be framed as the joint optimization of a suitably defined measure of the communication complexity of the distributed controller and its closed loop performance, in which these two competing metrics are traded off against each other in a principled manner.

In general one can select communication architectures that range in complexity from completely decentralized, i.e., distributed controllers with no communication allowed between sub-controllers, to essentially centralized and without delay, i.e., distributed controllers with instantaneous communication allowed between all sub-controllers.  However, if we ask that the distributed optimal controller restricted to the designed communication architecture be specified by the solution to a convex optimization problem then this limits the simplicity of the designed communication scheme \cite{W68,RL06,RCL10,LL11_QI}.  In particular a sufficient, and under mild assumptions necessary, condition for a distributed optimal controller to be specified by the solution to a convex optimization problem\footnote{For a more detailed overview of the relationship between information exchange constraints and the convexity of distributed optimal control problems, we refer the reader to \cite{RL06,RCL10,BV05,MMRY12} and the references therein.}  is that the communication architecture allow sub-controllers to communicate with each other as quickly as their control actions propagate through the plant \cite{RCL10}.  Although this condition may seem restrictive, it can often be met in practice by constructing a communication topology that mimics or is a superset of the physical topology of the plant.  For example, these delay based conditions \edit{may be satisfied in a smart-grid setting if fiber-optic cables are laid down in parallel to transmission lines; in a SDN setting if control packets are given priority in routing protocols; and in an automated highway system setting if vehicles are allowed to communicate wirelessly with nearby vehicles.}

When the aforementioned delay based condition is satisfied by a distributed constraint, it is said to 
be \emph{quadratically invariant} (QI) \cite{RL06,RCL10}.  \edit{While the resulting distributed optimal control problem is convex when quadratic invariance holds, it may still be infinite dimensional}.  Recently it has been shown that in the case of $\Htwo$ distributed optimal control subject to QI constraints imposed by a strongly connected communication architecture, i.e. one in which every sub-controller can exchange information with every other sub-controller subject to delay, the resulting distributed optimal controller synthesis problem can be reduced to a finite dimensional convex program, and hence admits an efficient solution \cite{LDXX,LD14}.\footnote{Other solutions exist to the $\Htwo$ distributed control problem subject to delay constraints -- we refer the reader to the discussion and references in \cite{LD14} for a more extensive overview of this literature.}  In light of these observations, we look to design \emph{strongly connected communication architectures} that induce \emph{QI} constraint sets -- once such a communication architecture is obtained, the methods from \cite{LDXX,LD14} can then be used to compute the optimal distributed controller restricted to that communication architecture exactly.

\edit{\textbf{Related prior work:}  Regularization techniques based on atomic norms have been employed to great success in system identification \cite{SBTR12,Fazel,MR_CDC14,LjungNewOld}.}
 \edit{As far as we are aware, the first instance of the use of regularization for the purpose of designing the architecture of a controller can be found in \cite{FLJ11} (these methods were then further developed in \cite{LFJ13}), in which an $\ell_1$ penalty is used with non-convex optimization to synthesize sparse static state feedback controllers with respect to an $\Htwo$ performance metric.  Other representative examples include the use of $\ell_1$ regularization to design  sparse treatment therapies \cite{JRM14}, consensus \cite{DLFM12,XB07} and synchronization \cite{FLJ13} topologies; and the use of group norm like penalties to design actuation/sensing schemes \cite{MC_CDC14,DJL14,Pol13}. }

\textbf{Contributions:}
\edit{We show that the communication complexity of a distributed controller can be inferred from the structure of its impulse response elements.  We use this observation to provide an explicit construction of an atomic norm \cite{CRPW12,Bon91,Pis86}, which we call the communication link norm, that can be incorporated into the RFD framework \cite{MC_RFD_TAC} to design strongly connected communication graphs that generate QI subspaces.}  As argued above, these two structural properties allow for the distributed optimal controller implemented using the designed communication architecture to be specified by the solution to a finite dimensional convex optimization problem \cite{LDXX,LD14}.  We also show that by augmenting the variational solution to the $\Htwo$ distributed optimal control problem presented in \cite{LDXX,LD14} with the communication link norm as a regularizer, the communication architecture/control law co-design problem can be formulated as a second order cone program.  By varying the weight on the communication link norm penalty function, the controller designer can use our co-design algorithm to explore the tradeoff between communication architecture complexity and closed loop performance in a principled way via convex optimization. We use these results to formulate a communication architecture/control law co-design algorithm that yields a distributed optimal controller and the communication architecture on which it is to be implemented. 

\textbf{Paper Organization:}
In \S \ref{sec:prelims} we introduce necessary operator theoretic concepts and establish notation.   In \S \ref{sec:problem} we formulate the communication architecture/control law co-design problem as the joint optimization of a suitably defined measure of the communication complexity of a distributed controller and the closed loop performance that it achieves.  In \S \ref{sec:delays}, we show how communication graphs can be used to generate distributed constraints, and show that if a communication graph that generates a QI subspace is augmented with additional communication links, the subspace generated by the resulting communication graph is also QI.  We use this observation and techniques from structured linear inverse problems \cite{CRPW12} in \S \ref{sec:co_design} to construct a convex regularizer that penalizes the use of additional communication links by a distributed controller, and formulate the co-design procedure.  In \S \ref{sec:examples} we discuss the computational complexity of the co-design procedure and illustrate the usefulness of our approach with two numerical examples.  We end with a discussion in \S \ref{sec:discussion}.

\section{Preliminaries}
\label{sec:prelims}
\subsection{Operator Theoretic Preliminaries}
\edit{We use standard definitions of the Hardy spaces $\Htwo$ and $\Hinf$.  We denote the restrictions of $\Hinf$ and $\Htwo$ to the space of real rational proper transfer matrices $\Rp$ by $\RHinf$ and $\RHtwo$, respectively.  As we work in discrete time, the two spaces are equal, and as a matter of convention  we refer to this space as $\RHinf$.  We refer the reader to \cite{ZDG96} for a review of this standard material.  For a signal $\vct f = (f^{(t)})_{t=0}^\infty$, we use $\vct f^{\leq d}$ to denote the truncation of $\vct f$ to its elements $f^{(t)}$ satisfying $t\leq d$, i.e., $\vct f^{\leq d} := (f^{(t)})_{t=0}^d$.  We extend the Banach space $\ell^{n}_2$ to the space
\begin{equation}
\ell^{n}_{2,e} := \{ \vct f:\Z_+ \to \R^{n} \, | \, \vct f^{\leq d} \in \ell^{n}_2 \text{ for all }  d\in \Z_+\},
\end{equation}
where $\Zp$ \edit{($\Zpp$)} denotes the set of non-negative \edit{(positive)} integers. A plant $G\in\Rp^{m\times n}$ can then be viewed as a linear map from $\ell^{n}_{2,e}$ to $\ell^{m}_{2,e}$.  Unless required, we do not explicitly denote dimensions and we assume that all vectors, operators and spaces are of compatible dimension throughout.}

\subsection{Notation}
We denote elements of $\ell_{2,e}$ with boldface lower case Latin letters, elements of $\Rp$ (which include matrices) with upper case Latin letters, and affine maps from $\RHinf$ to $\RHinf$ with upper case Fraktur letters such as $\mathfrak{M}$.  We denote
temporal indices, horizons and delays by lower case Latin letters. 



We denote the elements of the power series expansion of a map $G \in \RHinf$ by $G^{(t)}$, i.e., $G = \sum_{t=0}^\infty \frac{1}{z^t}G^{(t)}$.  \edit{We use $\RHinf^{\leq d}$ to denote the subspace of $\RHinf$ composed of finite impulse response (FIR) transfer matrices of horizon $d$, i.e., $\RHinf^{\leq d} := \{ G \in \RHinf \, {|} \, G = \sum_{t=0}^d \frac{1}{z^t}G^{(t)}\}$.  Similarly, we use $\RHinf^{\geq d+1}$ to denote the subspace of $\RHinf$ composed of transfer matrices with power series expansion elements satisfying $G^{(t)}=0$ for all $t\leq d$, i.e., $\RHinf^{\geq d+1} := \{ G \in \RHinf \, {|} \, G = \sum_{t=d+1}^\infty \frac{1}{z^t}G^{(t)}\}$.
For an element $G \in \RHinf$, we use $G^{\leq d}$ to denote the projection of $G$ onto $\RHinf^{\leq d}$, and $G^{\geq d+1}$ to denote the projection of $G$ onto $\RHinf^{\geq d+1}$, i.e., $G^{\leq d} = \sum_{t=0}^d \frac{1}{z^t}G^{(t)}$ and $G^{\geq d+1} = \sum_{t=d+1}^\infty \frac{1}{z^t}G^{(t)}$.  }

Sets are denoted by upper case script letters, such as $\mathscr{S}$, whereas subspaces of an inner product space are denoted by upper case calligraphic letters, such as $\s$.  
We denote the orthogonal complement of $\s$ with respect to the standard inner product on $\RHtwo$ by $\s^\perp$.  We use the greek letter $\Gamma$ to denote the adjacency matrix of a graph, and use labels in the subscript to distinguish among different graphs, i.e., $\Gamma_{\text{base}}$ and $\Gamma_{1}$ correspond to different graphs labeled ``base'' and ``1.''  We use $E_{ij}$ to denote the matrix with $(i,j)\text{th}$ element set to 1 and all others set to 0. We use $I_n$ and $0_n$ to denote the $n \times n$ dimensional identity matrix and all zeros matrix, respectively.  For a $p$ by $q$ block row by block column transfer matrix $M$ partitioned as $M=(M_{ij})$, we define the block support $\bsupp{M}$ of the transfer matrix $M$ to be the $p$ by $q$ integer matrix with $(i,j)$th element set to 1 if $M_{ij}$ is nonzero, and 0 otherwise.  Finally, we use the $\star$ superscript to denote that a parameter is the solution to an optimization problem.
%

%

\section{Communication Architecture Co-Design}
\label{sec:problem}
In this section we formulate the communication architecture/control law co-design problem as the joint optimization of a suitably defined measure of the communication complexity of the distributed controller and its closed loop performance. 
In particular, we introduce the convex optimization based solution to the $\Htwo$ distributed optimal control problem subject to delays presented in \cite{LDXX,LD14}, and modify this method to perform the communication architecture/control law co-design task.

\subsection{Distributed $\Htwo$ Optimal Control subject to Delays}
\label{sec:h2dist}
\begin{figure}[h!]
\centering
\includegraphics[width=1.5in]{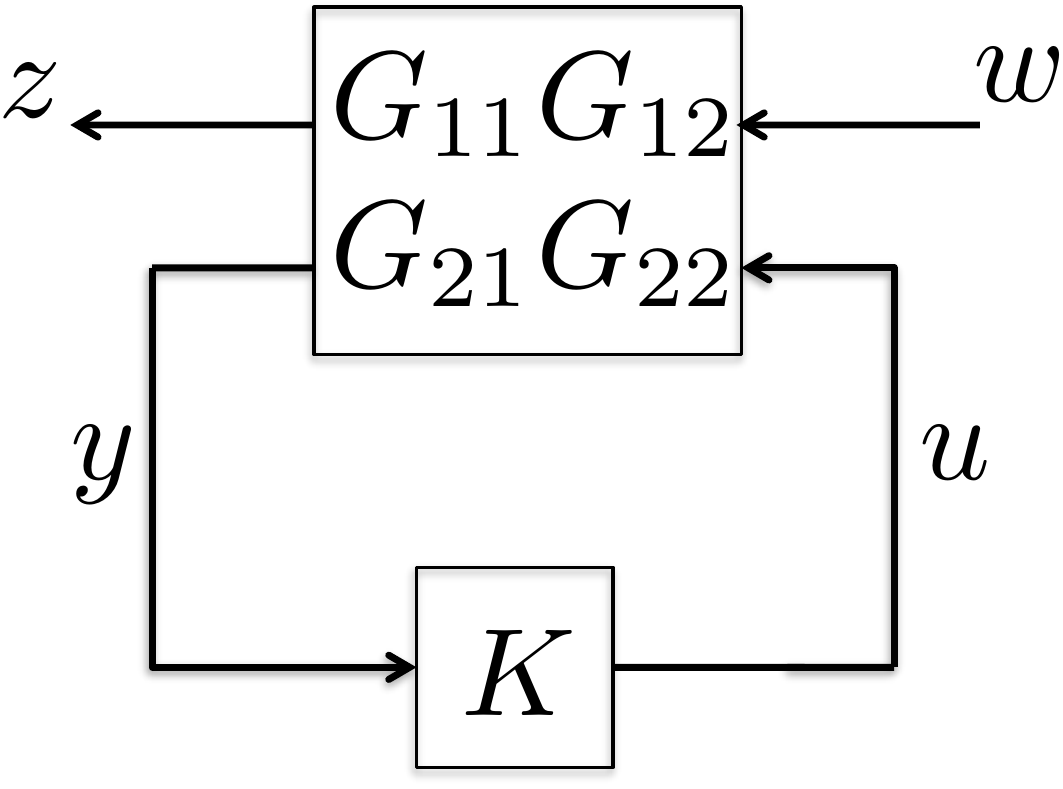}
\caption{A diagram of the generalized plant defined in \eqref{eq:gen_plant}.}
\label{fig:gen_plant}
\end{figure}

To review the relevant results of \cite{LDXX, LD14}, we introduce the discrete-time generalized plant $G$ given by 
\begin{equation}
G=\left[\begin{array}{c|cc}
A & B_{1} & B_{2}\\
\hline C_{1} & 0 & D_{12}\\
C_{2} & D_{21} & 0
\end{array}\right]=\begin{bmatrix}G_{11} & G_{12}\\
G_{21} & G_{22}
\end{bmatrix}
\label{eq:gen_plant}
\end{equation}
with inputs of dimension $p_{1},\, p_{2}$ and outputs of dimension
$q_{1},\, q_{2}$.  As illustrated in Figure \ref{fig:gen_plant}, this system describes the four transfer matrices from the disturbance and control inputs $\vct w$ and $\vct u$, respectively, to the controlled and measured outputs $\vct z$ and $\vct y$, respectively.
In order to ensure the existence of solutions to the necessary Riccati equations and to obtain simpler formulas, we assume that \edit{$(A,B_{1},C_1)$ and $(A,B_2,C_2)$ are both stabilizable and detectable, and that}
\begin{equation}
\edit{\text{ $D_{12}^{\top}D_{12}=I$, $D_{21}D_{21}^{\top}=I$, $C_{1}^{\top}D_{12}=0$,  $B_{1}D_{21}^{\top}=0$}.}
\label{eq:param_as}
\end{equation}  

Let $\s$ be a subspace that encodes the \edit{distributed constraints imposed on the controller $K$}.  For example, when some sub-controllers cannot access the measurements of other sub-controllers, the subspace $\s$ enforces corresponding sparsity constraints on the controller $K$. Alternatively, when sub-controllers can only gain access to other sub-controllers' measurements after a given delay, the subspace $\s$ enforces corresponding delay constraints on the controller $K$.

\edit{The distributed $\Htwo$ optimal control problem with subspace constraint $\s$ is then given by
\begin{equation}
\begin{array}{rl}
\minimize{K\in\Rp} & \Htwonorm{G_{11} - G_{12}K(I-G_{22}K)^{-1}G_{21}}^2 \\
\text{s.t.} & K \in \s\\
& \text{\edit{$K$ internally stabilizes $G$}}
\end{array}
\label{eq:distK_opt}
\end{equation}
where the objective function measures the $\Htwo$ norm of the closed loop transfer function from the exogenous disturbance $\vct w$ to the controlled output $\vct z$, and the first constraint ensures that the controller $K$ respects the distributed constraints imposed by the subspace $\s$.}

Optimization problem \eqref{eq:distK_opt} is in general both infinite dimensional and non-convex.  In \cite{LDXX, LD14}, the authors provide an exact and computationally tractable solution to optimization problem \eqref{eq:distK_opt} when the distributed constraint $\s$ is \emph{QI} \cite{RL06} with respect to $G_{22}$\footnote{A subspace $\s$ is said to be QI with respect to $G_{22}$ if $KG_{22}K \in \s$ for all $K \in \s$.  When quadratic invariance holds, we have that $K \in \s$ if and only if \edit{$K(I-G_{22}K)^{-1} \in \s$; }this key property allows for the convex parameterization \eqref{eq:dist_mm} of the distributed optimal control problem \eqref{eq:distK_opt}.} and is generated by a \emph{strongly connected} communication graph.  We say that a distributed constraint $\s$ is generated by a strongly connected communication graph\footnote{
We consider subspaces $\s$ that are strictly proper so that the reader can use the exact results presented in \cite{LD14}.  The authors of \cite{LD14} do however note that their method extends to non-strictly proper controllers at the expense of more complicated formulas.
} if it admits a decomposition of the form
\begin{equation}
\s = \Y \oplus \frac{1}{z^{d+1}}\Rp, \ \Y = \oplus_{t=1}^d \frac{1}{z^t}\Y^{(t)}
\label{eq:s_decomp}
\end{equation}
for some positive integer $d$, and some subspaces $\Y^{(t)} \subset \R^{p_2 \times q_2}$.  In \S\ref{sec:delays} we show how a strongly connected communication graph between sub-controllers can be used to define a subspace $\s$ that admits a decomposition \eqref{eq:s_decomp}. 

Restricting ourselves to distributed constraints $\s$ that are QI with respect to $G_{22}$ and that admit a decomposition of the form \eqref{eq:s_decomp} allows us to pose the optimal control problem \eqref{eq:distK_opt} as the following convex model matching problem
\edit{\begin{equation}
\begin{array}{rl}
\minimize{Q\edit{\in\RHinf}} & \Htwonorm{P_{11} - P_{12} Q P_{21}}^2 \\
\text{s.t.} & \mathfrak{C}\left(Q^{\leq d}\right) \in \Y \\
\end{array}
\label{eq:dist_mm}
\end{equation}}
through the use of a suitable Youla parameterization, where the $P_{ij}\in  \RHinf$ are appropriately defined stable transfer matrices and $\mathfrak{C}:\RHinf^{\leq d}\to\RHinf^{\leq d}$ is an appropriately defined affine map (cf.  \S III-B of \cite{LD14}).  It is further shown in \cite{LD14} that the solution $Q^\star$ to the distributed model matching problem \eqref{eq:dist_mm} with QI constraint $\s$ admitting decomposition \eqref{eq:s_decomp} is specified in terms of the solution to a finite dimensional convex quadratic program.  



\begin{thm}[Theorem 3 in \cite{LD14}]
\label{thm:quad}
Let $\s$ be QI under $G_{22}$ and admit a decomposition as in \eqref{eq:s_decomp}.   Let $Q^\star\in \s\cap\RHinf$ be the optimal solution to the convex model matching problem \eqref{eq:dist_mm}. 
Then $(Q^\star)^{\geq d+1}=0$ and
\begin{equation}
\begin{array}{rl}
(Q^\star)^{\leq d} = \underset{V\in\RHinf^{\leq d}}{\arg\min} & \Htwonorm{\mathfrak{L}\left(V\right)}^2 
\text{ s.t. }  \mathfrak{C}\left(V\right) \in \Y,
\end{array}
\label{eq:quad_prog}
\end{equation}
where $\mathfrak{L}$ is a linear map from $\RHinf^{\leq d}$ to $\RHinf^{\leq d}$, and $\mathfrak{C}$ is the affine map from $\RHinf^{\leq d}$ to $\RHinf^{\leq d}$ used to specify the model matching problem \eqref{eq:dist_mm}.
Furthermore, the optimal cost achieved by $Q^\star$ in the optimization problem \eqref{eq:dist_mm} is given by
\begin{equation}
\Htwonorm{P_{11}}^2+\Htwonorm{\mathfrak{L}\left((Q^\star)^{\leq d}\right)}^2.
\label{eq:cost}
\end{equation}
\end{thm}
\begin{rem}
The term $\Htwonorm{\mathfrak{L}\left((Q^\star)^{\leq d}\right)}^2$ in the optimal cost \eqref{eq:cost} quantifies the deviation of the performance achieved by the distributed optimal controller from that achieved by the centralized optimal controller.
\end{rem}

The optimization problem \eqref{eq:quad_prog} is finite dimensional because the maps $\mathfrak{L}$ and $\mathfrak{C}$ are both finite dimensional (they map the finite dimensional space $\RHinf^{\leq d}$ into itself) and act on the finite dimensional transfer matrix $V \in \RHinf^{\leq d}$.  These maps can be computed in terms of the state-space parameters of the generalized plant \eqref{eq:gen_plant} and the solution to appropriate Riccati equations (cf. \S III-B and \S IV-A of \cite{LD14}).  Under the assumptions \eqref{eq:param_as} the map $\mathfrak{L}$ is injective, and hence the convex quadratic program \eqref{eq:quad_prog} has a unique optimal solution $(Q^\star)^{\leq d}$.

As the distributed constraint $\s$ is assumed to be QI, the optimal distributed controller $K^\star \in \s$ specified by the solution to the non-convex optimization problem \eqref{eq:distK_opt} can be recovered from the optimal Youla parameter $Q^\star\in \s$ through a suitable linear fractional transformation (cf. Theorem 3 of \cite{LD14}). 

\begin{rem}
If the state-space matrix $A$ specified in the generalized plant \eqref{eq:gen_plant} is of dimension $s \times s$, then the resulting optimal controller $K^\star$ admits a state-space realization of order $s+q_2d$.  As argued in \cite{LD14}, this is at worst within a constant factor of the minimal realization order.
\end{rem}

\subsection{Communication \edit{Delay Co-Design} via Convex Optimization}
Although our objective is to design the communication graph on which the distributed controller $K$ is implemented, for the computational reasons described in \S\ref{sec:h2dist} it is preferable to solve a problem in terms of the Youla parameter $Q$ as this leads to the convex optimization problems \eqref{eq:dist_mm} and \eqref{eq:quad_prog}. 
In order to perform the communication architecture/control law co-design task in the Youla domain, we restrict ourselves to designing strongly connected communication architectures that generate QI subspaces, i.e., subspaces that are QI and that admit a decomposition of the form \eqref{eq:s_decomp}.  As argued in \S\ref{sec:introduction}, this is a practically relevant class of communication architectures to consider, and further, based on the previous discussion it is then possible to solve for the resulting distributed optimal controller restricted to the designed communication architecture using the results of Theorem \ref{thm:quad}.  

Our approach to accomplish the co-design task is to remove the subspace constraint $\mathfrak{C}\left(V\right) \in \Y$, which encodes the distributed structure of the controller, from the optimization problem \eqref{eq:quad_prog} and to augment the objective of the optimization problem with a convex penalty function that instead induces suitable structure in $\mathfrak{C}\left(V\right)$.  
In particular, we seek a convex penalty function $\commnorm{\cdot}$ \edit{and horizon $d$} such that the structure of $\mathfrak{C}\left(V^\star\right)$, where $V^\star$ is the solution to 
\edit{\begin{equation}
\begin{array}{rl}
\minimize{V\in\RHinf^{\leq d}} & \Htwonorm{\mathfrak{L}\left(V\right)}^2 + \lambda \commnorm{\mathfrak{C}\left(V\right)},
\end{array}
\label{eq:reg_opt}
\end{equation}}
\edit{can be used to define} an appropriate QI subspace $\s$ that admits a decomposition of the form \eqref{eq:s_decomp}.  Imposing that the designed subspace $\s$ be QI ensures that the structure induced in $\mathfrak{C}\left(V^\star\right)$ corresponds to the structure of the resulting distributed controller $K^\star$.  Further imposing that the designed subspace $\s$ admit a decomposition of the form \eqref{eq:s_decomp} ensures that the distributed optimal controller restricted to lie in the subspace $\s$ can be computed using Theorem \ref{thm:quad}.
\begin{rem}
The regularization weight $\lambda\geq 0$ allows the controller designer to tradeoff between closed loop performance (as measured by $\Htwonorm{\mathfrak{L}\left(V\right)}^2$) and communication complexity (as measured by $\commnorm{\mathfrak{C}\left(V\right)}$).  
\end{rem}

In order to define an appropriate convex penalty $\commnorm{\cdot}$, we need to understand how a communication graph between sub-controllers defines the subspace $\Y$ in which $\mathfrak{C}\left(V\right)$  is constrained to lie in optimization problem \eqref{eq:quad_prog} -- this in turn informs what structure to induce in $\mathfrak{C}\left(V^\star\right)$ in the regularized optimization problem \eqref{eq:reg_opt}.  To that end, in \S \ref{sec:delays} we define a simple communication protocol between sub-controllers that allows communication graphs to be associated with distributed subspace constraints in a natural way.  Within this framework, we show that if a communication graph generates a distributed subspace $\s$ that is QI with respect to $G_{22}$, then adding additional communication links to this graph preserves the QI property of the distributed subspace that it generates.  We use this observation to pose the communication architecture design problem as one of augmenting a suitably defined base communication graph, namely a simple graph that generates a QI subspace, with additional communication links.


\section{Communication Graphs and Quadratically Invariant Subspaces}
\label{sec:delays}
This section first shows how a communication graph connecting sub-controllers can be used to define the subspace $\s$ in which the controller $K$ is constrained to lie in the distributed optimal control problem \eqref{eq:distK_opt}.  In particular, if two sub-controllers exchange information using the shortest path between them on an underlying communication graph, then there is a natural way of generating a subspace constraint from the adjacency matrix of that graph.  Under this information exchange protocol, we then define a set of strongly connected communication graphs that generate subspace constraints that are QI with respect to a plant $G_{22}$ in terms of a \emph{base} and a \emph{maximal} communication graph.  This approach allows the controller designer to specify which communication links between sub-controllers are \emph{physically realizable}, i.e., which communication links can be built subject to the physical constraints of the system.

\subsection{Generating Subspaces from Communication Graphs}
\label{sec:graphs}
Consider a generalized plant \eqref{eq:gen_plant} comprised of $n$ sub-plants, each equipped with its own sub-controller.  Let $\N:=\{1,\dots,n\}$ and label each sub-controller by a number $i \in \N$.  To each such sub-controller $i$ associate a space of possible control actions $\mathcal{U}_i=\ell^{p_{2,i}}_{2,e}$ and a space of possible output measurements $\mathcal{Y}_i = \ell^{q_{2,i}}_{2,e}$, and define the overall control and measurement spaces as  $\mathcal{U}:=\mathcal{U}_1 \times \cdots \times \mathcal{U}_n$ and $\mathcal{Y}:=\mathcal{Y}_1 \times \cdots \times \mathcal{Y}_n$, respectively.

Then, for any pair of sub-controllers $i$ and $j$, the $(i,j)^\text{th}$ block of $G_{22}$ is the mapping from the control action $\vct u_j$ taken by sub-controller $j$ to the measurement $\vct y_i$ of sub-controller $i$, i.e., $\left(G_{22}\right)_{ij}:\mathcal{U}_j\to\mathcal{Y}_i$.  Similarly, the mapping from the measurement $\vct y_j$, transmitted by sub-controller $j$, to the control action $\vct u_i$ taken by sub-controller $i$ is given by $K_{ij}:\mathcal{Y}_j \to \mathcal{U}_i$.

We then form the overall measurement and control vectors
\begin{equation}
\vct y = \begin{bmatrix} (\vct y_1)^\top & \cdots & (\vct y_n)^\top \end{bmatrix}^\top, \indent \vct u = \begin{bmatrix} (\vct u_1)^\top & \cdots & (\vct u_n)^\top \end{bmatrix}^\top
\label{eq:iOS}
\end{equation}
leading to the natural block-wise partitions of the plant $G_{22}$
\begin{equation}
G_{22} = \begin{bmatrix} \left(G_{22}\right)_{11} & \cdots &  \left(G_{22}\right)_{1n} \\
\vdots & \ddots & \vdots \\
 \left(G_{22}\right)_{n1} & \cdots & \left(G_{22}\right)_{nn}\end{bmatrix}
 \label{eq:splitG}
 \end{equation}
 and of the controller $K$
\begin{equation}
K = \begin{bmatrix} K_{11} & \cdots & K_{1n} \\
\vdots & \ddots & \vdots \\
K_{n1} & \cdots & K_{nn}\end{bmatrix}.
\label{eq:splitK}
\end{equation}

We assume that sub-controllers exchange measurements with each other subject to delays imposed by an underlying communication graph -- specifically, we assume that sub-controller $i$ has access to sub-controller $j$'s measurement $\vct y_j$ with delay specified by the length of the shortest path from sub-controller $j$ to sub-controller $i$ in the communication graph.  Formally, let $\Gamma$ be the adjacency matrix of the communication graph between sub-controllers, i.e., $\Gamma$ is the integer matrix with rows and columns indexed by $\N$, such that $\Gamma_{kl}$ is equal to 1 if there is an edge from $l$ to $k$, and 0 otherwise.  The \emph{communication delay} from sub-controller $j$ to sub-controller $i$ is then given by the length of the shortest path from $j$ to $i$ as specified by the adjacency matrix gamma $\Gamma$. In particular, we define\footnote{See Lemma 8.1.2 of \cite{GR01_AGT} for a graph theoretic justification of this definition.} the communication delay from sub-controller $j$ to sub-controller $i$ to be given by
\begin{equation} 
\edit{c_{ij} := \min\left\{ d \in \Zp \, \big{|} \, \Gamma_{ij}^d \neq 0\right\}}
\label{eq:comm_delays}
\end{equation}
\edit{if an integer satisfying the condition in \eqref{eq:comm_delays} exists, and set $c_{ij} = \infty$ otherwise.}

We say that a strictly proper distributed controller $K$ can be \emph{implemented} on a communication graph with adjacency matrix $\Gamma$ if for all $i,j\in \N$, we have that the the $(i,j)$th block of the controller $K$ satisfies $K_{ij}^{(t)} = 0$ for all positive integers $t \leq c_{ij}$, or equivalently, that $K_{ij} \in \frac{1}{z^{c_{ij}+1}}\Rp$.  In words, this says that sub-controller $j$ only has access to the measurement $\vct y_i$ from sub-controller $i$ after $c_{ij}$ time steps, the length of the shortest path from $j$ to $i$ in the communication graph, and can only take actions based on this measurement after a computational delay of one time step.\footnote{This computational delay is included to ensure that the resulting controller is strictly proper.}
More succinctly, this condition holds if 
$\bsupp{K^{(t)}} \subseteq \supp{\Gamma^{t-1}}$ for all $t \geq 1$.

If $\Gamma$ is the adjacency matrix of a strongly connected graph, then there exists a path between all ordered pairs of sub-controllers $(i,j) \in \N\times\N$ --  this implies that there exists a positive delay $d(\Gamma)$ after which a given measurement $\vct y_j$ is available to all sub-controllers.  
  In particular, we define the delay $d(\Gamma)$ associated with the adjacency matrix $\Gamma$ to be
\begin{equation}
\edit{d\left(\Gamma\right) := \sup\left\{ \tau\in\Zpp \, \big{|} \, \exists (k,l) \in \N\times\N \text{ s.t. } \Gamma^{\tau-1}_{kl} = 0 \right\}.}
\label{eq:delta}
\end{equation}
Using this convention all measurements $y_j^{(t)}$ are available to all sub-controllers by time $t+d(\Gamma)+1$.  When the delay $d(\Gamma)$ is finite, we say that $\Gamma$ is a \emph{strongly connected adjacency matrix}, as it defines a strongly connected communication graph.

We define the subspace $\s(\Gamma)$ generated by a strongly connected adjacency matrix $\Gamma$ to be
\begin{equation}
\s(\Gamma) := \Y(\Gamma)\oplus \frac{1}{z^{d(\Gamma)+1}}\Rp,
\label{eq:s_gam}
\end{equation}
where $d(\Gamma)$ is as defined in \eqref{eq:delta}, and
$\Y(\Gamma) := \oplus_{t=1}^d \frac{1}{z^t} \Y^{(t)}(\Gamma)$ is specified by the subspaces
\begin{equation}
\Y^{(t)}(\Gamma) := \left\{ M \in \R^{p_2 \times q_2} \, \big{|} \, \bsupp{M} \subseteq \supp{\Gamma^{t-1}}\right\}.
\label{eq:Yts}
\end{equation}
It is then immediate that a controller $K$ can be implemented on the communication graph $\Gamma$ if and only if $K \in \s(\Gamma)$.

\begin{figure}[t]
\centering
\includegraphics[width = .15\textwidth]{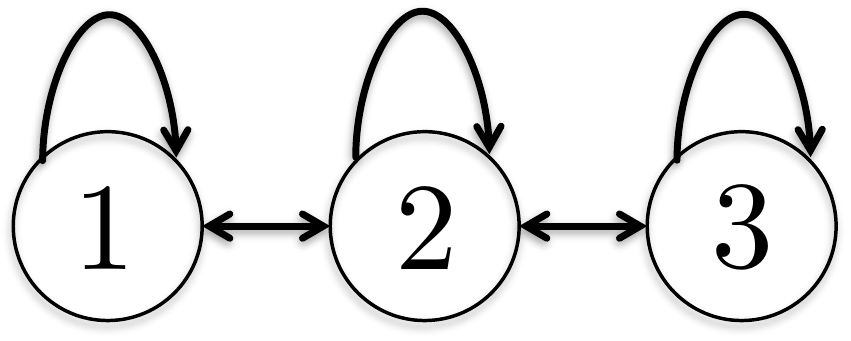}
\caption{Three subsystem chain example}
\label{fig:4chain}
\end{figure}
\begin{example}
Consider the communication graph illustrated in Figure \ref{fig:4chain} with strongly connected adjacency matrix $\Gamma_{\text{3-chain}}$ given by 
\begin{equation}
\Gamma_{\text{3-chain}} = \begin{bmatrix} 1 & 1 & 0 \\ 1 & 1 & 1 \\ 0 & 1 & 1 \end{bmatrix}. \label{eq:3-chain_gam}
\end{equation}
This communication graph generates the subspace
\begin{equation}
\s(\Gamma_{\text{3-chain}}):= \frac{1}{z}\begin{bmatrix} \ast & 0 & 0  \\ 0 & \ast & 0 \\ 0 & 0 & \ast \end{bmatrix} \oplus 
\frac{1}{z^2}\begin{bmatrix} \ast & \ast & 0  \\  \ast & \ast & \ast \\ 0 & \ast & \ast  \end{bmatrix} \oplus \frac{1}{z^3}\Rp,
\label{eq:3-chain_s}
\end{equation}
where $\ast$ is used to denote a space of appropriately sized real matrices.  The communication delays associated with this graph are then given by $c_{ij} = |i-j|$ (e.g., $c_{11} = 0$, $c_{12} = 1$ and $c_{13} = 2$).  We also have that $d(\Gamma_{\text{3-chain}})= 2$, which is the length of the longest path between nodes in this graph, and that
\[
\Y(\Gamma_{\text{3-chain}})= \frac{1}{z}\begin{bmatrix} \ast & 0 & 0  \\ 0 & \ast & 0 \\ 0 & 0 & \ast \end{bmatrix} \oplus 
\frac{1}{z^2}\begin{bmatrix} \ast & \ast & 0  \\  \ast & \ast & \ast \\ 0 & \ast & \ast  \end{bmatrix}\subset \RHinf^{\leq 2}.
\]
\hfill{}\QED
\label{ex:3chain_s_and_y}
\end{example}

Thus, given such a strongly connected adjacency matrix  $\Gamma$, the distributed optimal controller $K^\star$ implemented using the graph specified by $\Gamma$ can be obtained by solving the optimization problem \eqref{eq:distK_opt} with subspace constraint $\s(\Gamma)$  -- however, this optimization problem can only be reformulated as the convex programs \eqref{eq:dist_mm} and \eqref{eq:quad_prog} if the subspace $\s(\Gamma)$ is QI with respect to $G_{22}$ \cite{LL11_QI}. 

\subsection{Quadratically Invariant Communication Graphs}
The discussion of \S\ref{sec:problem} and \S\ref{sec:graphs} shows that communication graphs that are strongly connected and that generate a subspace \eqref{eq:s_gam} that is QI with respect to $G_{22}$ \edit{allow for the distributed optimal control problem \eqref{eq:distK_opt} to be solved via the finite dimensional convex program \eqref{eq:quad_prog}.}  In this subsection, we  characterize a set of such communication graphs in terms of a \emph{base QI} and a \emph{maximal QI communication graph} corresponding to a plant $G_{22}$.  The base QI communication graph defines a simple communication architecture that generates a QI subspace, whereas the maximal QI communication graph is the densest communication architecture that can be built given the physical constraints of the system. 

We assume that the sub-controllers have disjoint measurement and actuation channels, i.e., that $B_2$ and $C_2$ are block-diagonal, and that the dynamics of the system are strongly connected, i.e., that $\bsupp{A}$ corresponds to the adjacency matrix of a strongly connected graph.  We discuss alternative approaches for when these assumptions do not hold in \S\ref{sec:discussion}.  For the sake of brevity, we often refer to a communication graph by its adjacency matrix $\Gamma$.


\subsubsection*{The base QI communication graph}  Our objective is to identify a simple communication graph, i.e., a graph defined by a sparse adjacency matrix $\Gamma_{\text{base}}$, such that the resulting subspace $\s(\Gamma_{\text{base}})$ is QI with respect to $G_{22}$.  To that end,  let the \emph{base QI communication graph} of plant $G_{22}$ with realization \eqref{eq:gen_plant} be specified by the adjacency matrix
\begin{equation}
\Gamma_{\text{base}}:=\bsupp{A}.
\label{eq:gam_min}
\end{equation}
Notice that under the block-diagonal assumptions imposed on the state-space parameters $B_{2}$ and $C_{2}$, this implies that $\Gamma_{\text{base}}$ mimics or is a superset of the physical topology of the plant $G_{22}$, as \edit{$\bsupp{G_{22}^{(t)}} = \bsupp{C_2 A^{t-1} B_2} \subseteq \bsupp{A}^{t-1}$.}

Define the \emph{propagation delay} from sub-plant $j$ to sub-plant $i$ of a plant $G_{22}$ to be the largest integer $p_{ij}$ such that 
\begin{equation}
\left(G_{22}\right)_{ij} \in \frac{1}{z^{p_{ij}}}\Rp.
\label{eq:prop_delays}
\end{equation}
It is shown in \cite{RCL10} that if a subspace $\s$ constrains the blocks of the controller $K$ to satisfy $K_{kl} \in \frac{1}{z^{c_{kl}+1}}\Rp$, and the communication delays\footnote{These are equivalent to the prior definition \eqref{eq:comm_delays} of communication delays $\{c_{kl}\}$.} $\{c_{kl}\}$ satisfy the triangle inequality $c_{ki} +c_{ij} \geq c_{kj}$, then $\s$ is QI with respect to $G_{22}$ if
\begin{equation}
c_{ij} \leq p_{ij} +1
\label{eq:delay_cdt}
\end{equation} for all $i,j \in \N$.  An intuitive interpretation of this condition is that $\s$ is QI if it allows sub-controllers to communicate with each other as fast as their control actions propagate through the plant.  Since we take the base QI communication graph $\Gamma_{\text{base}}$ to mimic the topology of the plant $G_{22}$, we expect this condition to hold and for $\s(\Gamma_{\text{base}})$ to be QI with respect to $G_{22}$.  We formalize this intuition in the following lemma.
\begin{lem}
Let the plant $G_{22}$ be specified by state-space parameters $(A,B_2,C_2)$, and suppose that $B_2$ and $C_2$ are block diagonal.  Let $\{p_{ij}\}$ denote the propagation delays of the plant $G_{22}$ as defined in \eqref{eq:prop_delays}.  Assume that $\Gamma_{\text{base}}$, as specified as in equation \eqref{eq:gam_min}, is a strongly connected adjacency matrix, and let $\{b_{ij}\}$ denote the communication delays \eqref{eq:comm_delays} imposed by the adjacency matrix $\Gamma_{\text{base}}$.  The communication delays $\{b_{ij}\}$ then satisfy condition \eqref{eq:delay_cdt} and the subspace $\s(\Gamma_{\text{base}})$ is quadratically invariant with respect to $G_{22}$.
\label{lem:base_delay}
\end{lem}
\begin{IEEEproof}
The definition of the base QI communication graph $\Gamma_{\text{base}}$ and the assumption that $B_2$ and $C_2$ are block-diagonal imply that $\bsupp{G_{22}^{(t)}} \subseteq \bsupp{A^{t-1}}\subseteq\supp{\Gamma_{\text{base}}^{t-1}}$.  This in turn can be verified to guarantee that \eqref{eq:delay_cdt} holds.  \edit{Thus it suffices to show that the communication delays $\{b_{kl}\}$ satisfy the triangle inequality $b_{ki} +b_{ij} \geq b_{kj}$ for all $i,j,k \in \N$. 
First observe that (i) $b_{ii} + b_{ii} \geq b_{ii}$, and (ii) $b_{ii} + b_{ij} \geq b_{ij}$, as all $b_{ij}\geq 0$.  Thus it remains to show that $b_{ki} + b_{ij} \geq b_{kj}$ for $i\neq j \neq k$.  Suppose, seeking contradiction, that 
\begin{equation}
b_{ki} + b_{ij} < b_{kj}.
\label{eq:triangle1}
\end{equation}}
Note that by definition \eqref{eq:comm_delays} of the communication delays and Lemma 8.1.2 of \cite{GR01_AGT}, the  inequality \eqref{eq:triangle1} is equivalent to
\begin{equation}
\begin{array}{l}
\min\{r \, | \, \exists\text{ path of length $r$ from $i$ to $k$}\} + \\ \indent \min\{r \, | \, \exists\text{ path of length $r$ from $j$ to $i$}\}  <\\  \indent \indent  \min\{r \, | \, \exists\text{ path of length $r$ from $j$ to $k$}\}.
\end{array}
\end{equation}
Notice however that we must have that
\begin{equation}
\begin{array}{l}\min\{r \, | \, \exists\text{ path of length $r$ from $j$ to $k$}\} \leq \\ \indent \min\{r \, | \, \exists\text{ path of length $r$ from $j$ to $i$}\} + \\ \indent \indent \min\{r \, | \, \exists\text{ path of length $r$ from $i$ to $k$}\},
\end{array}
\label{eq:triangle2}
\end{equation}
as the concatenation of a path from $j$ to $i$ and a path from $i$ to $k$ yields a path from $j$ to $k$.  Combining inequalities \eqref{eq:triangle1} and \eqref{eq:triangle2} yields the desired contradiction, proving the result.
\end{IEEEproof}

Lemma \ref{lem:base_delay} thus provides a simple means of constructing a base QI communication graph by taking a communication topology that mimics the physical topology of the plant $G_{22}$.

\subsubsection*{Augmenting the base QI communication graph}  The delay condition \eqref{eq:delay_cdt} suggests that a natural way of constructing QI communication architectures given a base QI communication graph is to augment the base graph with additional communication links, as adding a link to a communication graph can only decrease its communication delays $c_{ij}$.  
\begin{prop}
Let $\Gamma_{\text{base}}$ be defined as in \eqref{eq:gam_min}, and let $\Gamma$ be an adjacency matrix satisfying $\supp{\Gamma_{\text{base}}}\subset\supp{\Gamma}$.  Then the generated subspace $\s(\Gamma)$, as defined in \eqref{eq:s_gam}, is quadratically invariant with respect to $G_{22}$.
\label{prop:augment}
\end{prop}
\begin{IEEEproof}
Let $\{b_{ij}\}$ and $\{c_{ij}\}$ denote the communication delays associated with the base QI communication graph $\Gamma_{\text{base}}$ and the augmented communication graph $\Gamma$, respectively.  It follows from the definition of the communication delays $\eqref{eq:comm_delays}$ that the support nesting condition $\supp{\Gamma_{\text{base}}}\subset\supp{\Gamma}$ implies that $b_{ij}\geq c_{ij}$ for all $i,j\in\N$.   By Lemma \ref{lem:base_delay} we have that $b_{ij}\leq p_{ij}+1$, and therefore $c_{ij} \leq b_{ij} \leq p_{ij}+1$.  An identical argument to that used to prove Lemma \ref{lem:base_delay} shows that the delays $c_{ij}$ satisfy the required triangle inequality, implying that $\s(\Gamma)$ is QI with respect to $G_{22}$.
\end{IEEEproof}

In words, the nesting condition  $\supp{\Gamma_{\text{base}}}\subset\supp{\Gamma}$ simply means that the communication graph $\Gamma$ can be constructed by adding communication links to the base QI communication graph $\Gamma_{\text{base}}$.  It follows that any graph built by augmenting $\Gamma_{\text{base}}$ with additional communication links generates a QI subspace \eqref{eq:s_gam}.

\edit{
\begin{rem}
Although we have suggested a specific construction for $\Gamma_\text{base}$, Proposition \ref{prop:augment} makes clear that any strongly connected graph that generates a subspace constraint that is QI with respect to $G_{22}$ can be used as the base QI communication graph.  We discuss the implications of this added flexibility in \S\ref{sec:discussion}.
\label{rem:base}
\end{rem}}

\subsubsection*{The maximal QI communication graph}  In order to augment the base QI communication graph in a physically relevant way, one must first specify what additional communication links can be built given the physical constraints of the system.  For example, if two sub-controllers are separated by a large physical distance, it may not be possible to build a direct communication link between them.  The set of \edit{additional} communication links that can be physically constructed is application dependent -- we therefore assume that the controller designer has specified a collection $\mathscr{E}$ of directed edges that define what communication links can be built \edit{in addition to those already present in the base QI communication graph}.  \edit{In particular, we assume that it is possible to build a direct communication link from sub-controller $j$ to sub-controller $i$, i.e., to build a communication graph $\Gamma_{\text{built}}=\Gamma_\text{base}+\Gamma$ with $\Gamma_{ij}=1$, only if $(i,j) \in \mathscr{E}$.} 

Given a collection of directed edges $\mathscr{E}$, the \emph{maximal QI communication graph} $\Gamma_{\max}$ is given by
\begin{equation}
\Gamma_{\max} := \Gamma_{\text{base}} + M,
\label{eq:gam_max}
\end{equation} 
where $M$ is a $n \times n$ dimensional matrix with $M_{ij}$ set to 1 if $(i,j) \in \mathscr{E}$ and 0 otherwise.  In words, the maximal QI adjacency matrix $\Gamma_{\max}$ specifies a communication graph that uses all possible communication links listed in the set $\mathscr{E}$, in addition to those links already used by the base QI communication graph.
Consequently, we say that a communication graph can be \emph{physically built} if its adjacency matrix $\Gamma$ satisfies
\begin{equation}
\supp{\Gamma} \subseteq \supp{\Gamma_{\max}},
\label{eq:phys_real}
\end{equation}
i.e., if it can be built from communication links used by the base QI communication graph and/or those listed in the set $\mathscr{E}$.

\subsubsection*{The QI communication graph design set}   We now define a set of strongly connected and physically realizable communication graphs that generate QI subspace constraints as specified in equation \eqref{eq:s_gam} -- in particular, the base and maximal QI graphs correspond to the boundary points of this set.

\begin{prop}
Given a plant $G_{22}$ and a set of directed edges $\mathscr{E}$, let the adjacency matrices $\Gamma_{\text{base}}$ and $\Gamma_{\max}$ of the base and maximal QI communication graphs be defined as in \eqref{eq:gam_min} and \eqref{eq:gam_max}, respectively.  Then an adjacency matrix $\Gamma$ corresponds to a strongly connected communication graph that can be physically built and that generates a quadratically invariant subspace $\s(\Gamma)$ of the form \eqref{eq:s_gam} if
\begin{equation}
\supp{\Gamma_{\text{base}}}\subseteq\supp{\Gamma}\subseteq\supp{\Gamma_{\max}}.
\label{eq:gam_space}
\end{equation}
\label{prop:space}
\end{prop}
\begin{IEEEproof}
Follows from Prop. \ref{prop:augment} and definitions \eqref{eq:gam_max} and \eqref{eq:phys_real}.
\end{IEEEproof}
The following corollary is then immediate.

\begin{coro}
Let $\Gamma_1$ and $\Gamma_2$ be adjacency matrices that satisfy the nesting condition \eqref{eq:gam_space} and suppose further that $\supp{\Gamma_1}\subseteq\supp{\Gamma_2}$.  Let $\nu_{\bullet}$, with $\bullet\in \left\{{\text{base}},{1}, {2}, {\max}\right\}$ be the closed loop norm achieved by the optimal distributed controller implemented using communication graph $\Gamma_{\bullet}$.
Then\edit{
\begin{equation}
d(\Gamma_{\text{base}}) \geq d(\Gamma_1) \geq d(\Gamma_2) \geq d(\Gamma_{\max}),
\label{eq:delta_nest}
\end{equation}}
\vspace{-4mm}
\begin{equation}
\s(\Gamma_{\text{base}}) \subseteq \s(\Gamma_1) \subseteq \s(\Gamma_2) \subseteq \s(\Gamma_{\max}),
\label{eq:s_nest}
\end{equation} and\edit{
\begin{equation} 
\nu_{\text{base}}\geq\nu_{1}\geq\nu_{2}\geq\nu_{\max}
\label{eq:norm_nest}
\end{equation}}
\label{coro:nest}
\end{coro}
\begin{IEEEproof}
Relations \eqref{eq:delta_nest} and \eqref{eq:s_nest} follow immediately from the hypotheses of the corollary and the definitions of the delays $d(\Gamma_\bullet)$ and the subspaces $\s(\Gamma_\bullet)$ as given in \eqref{eq:delta} and \eqref{eq:s_gam}, respectively.  The condition \eqref{eq:norm_nest} on the norms $\nu_\bullet$ follows immediately from the subspace nesting condition \eqref{eq:s_nest} and the fact that the optimal norm $\nu_\bullet$ achievable by a distributed controller implemented using a communication graph with adjacency matrix $\Gamma_\bullet$ is specified by the optimal value of the objective function of the optimization problem \eqref{eq:distK_opt} with distributed constraint $\s(\Gamma_\bullet)$.
\end{IEEEproof}

Corollary \ref{coro:nest} states that as more edges are added to the base QI communication graph, the performance of the optimal distributed controller implemented on the resulting communication graph improves.  Thus there is a quantifiable tradeoff between the communication complexity and the closed loop performance of the resulting distributed optimal controller.
%
 To fully explore this tradeoff, the controller designer would have to enumerate the \emph{QI communication graph design set} which is composed of adjacency matrices satisfying the nesting condition \eqref{eq:gam_space}.  Denoting this set by $\mathscr{G}$, a simple computation shows that \edit{$|\mathscr{G}| = 2^{|\mathscr{E}|}$} -- \edit{ thus the controller designer has to consider a set of graphs of cardinality exponential in the number of possible additional communication links}.  This poor scaling motivates the need for a principled approach to exploring the design space of communication graphs via the regularized optimization problem \eqref{eq:reg_opt}.

\section{The Communication Graph Co-Design Algorithm}
\label{sec:co_design}
In this section we leverage Propositions \ref{prop:augment} and \ref{prop:space} as well as tools from approximation theory \cite{CRPW12}, \approxcitations$\text{ }$  to construct a convex penalty function $\commnorm{\cdot}$, which we call the \emph{communication link norm}, that allows the controller designer to explore the QI communication graph design set $\mathscr{G}$ in a principled manner via the regularized convex optimization problem \eqref{eq:reg_opt}.  We then propose a communication architecture/control law co-design algorithm based on this optimization problem and show that it indeed does produce strongly connected communication graphs that generate quadratically invariant subspaces.
\subsection{The Communication Link Norm}
\label{sec:comm_norm}\edit{
Recall that our approach to the co-design task is to induce suitable structure in the expression $\mathfrak{C}\left(V^\star\right)$, where $V^\star$ is the solution to the regularized convex optimization problem \eqref{eq:reg_opt} employing the yet to be specified convex penalty function $\commnorm{\cdot}$.  We argued that the structure induced in the expression $\mathfrak{C}\left(V^\star\right)$ should correspond to a strongly connected communication graph that generates a QI subspace of the form \eqref{eq:s_decomp}, and characterized a set of graphs satisfying these properties, namely the QI communication graph design set $\mathscr{G}$.
To explore the QI communication graph design set $\mathscr{G}$, we begin with the base QI communication graph $\Gamma_{\text{base}}$ and augment it with additional communication links drawn from the set $\mathscr{E}$.  The convex penalty function $\commnorm{\cdot}$ used in the regularized optimization problem \eqref{eq:reg_opt} should therefore penalize the use of such additional communication links -- in this way the controller designer can tradeoff between communication complexity and closed loop performance by varying
the regularization weight $\lambda$ in optimization problem \eqref{eq:reg_opt}.}

We view distributed controllers implemented using a dense communication graph as being composed of a superposition of simple \emph{atomic} controllers that are implemented using simple communication graphs, i.e., using communication graphs obtained by adding a small number of edges to the base QI communication graph.  This viewpoint suggests choosing the convex penalty function $\commnorm{\cdot}$ to be an atomic norm \cite{CRPW12,Bon91,Pis86}.

Indeed, if one seeks a solution $X^\star$ that can be composed as a linear combination of a small number of atoms drawn from a set $\A$, then a useful approach, as described in  \cite{CRPW12,CDS98,CRT06,Donoho04, fazelThesis,RFP10,CR12}, to induce such structure in the solution of an optimization problem is to employ a convex penalty function that is given by the atomic norm induced by the atoms $\A$ \cite{Bon91,Pis86}. Examples of the types of structured solutions one may desire include sparse, group sparse and signed vectors,  and low-rank, permutation and orthogonal matrices \cite{CRPW12}.  Specifically, if one desires a solution $X^\star$ 
that admits a decomposition of the form
\begin{equation}
X^\star=\sum_{i=1}^{r}c_{i}A_{i},\ A_{i}\in\A,\ c_{i}\geq0
\label{eq:decomp}
\end{equation}
for a set of appropriately scaled and centered atoms $\A$,
and a small number $r$ relative to the ambient dimension, then solving
\begin{equation}
\begin{array}{rl}
\underset{X}{\mathrm{minimize}}  &  \Htwonorm{\mathfrak{A}(X)}^2 + \lambda\|X\|_\A
\end{array}
\end{equation}
with $\mathfrak{A}(\cdot)$ an affine map, and the atomic norm $\|\cdot\|_\A$ \edit{given by\footnote{\edit{If no such decomposition exists, then $\|X\|_\A = \infty$.}} 
\begin{equation}\begin{array}{c}
\|X\|_{\A}: = \inf\left\{\sum_{A\in\A}c_{A}\, \big{|}\, X=\sum_{A\in\A}c_{A}A,\, c_A \geq 0\right\}\end{array}
\label{eq:atom_norm}
\end{equation}}
results in solutions that are both consistent with the data as measured in terms of the cost function $\Htwonorm{\mathfrak{A}(X)}^2$, and that \edit{ admit sparse atomic decompositions, i.e., that are a combination of a small number of elements from $\A$. }

We can therefore fully characterize our desired convex penalty function $\commnorm{\cdot}$ by specifying its defining atomic set $\A_{\text{comm}}$ and then invoking definition \eqref{eq:atom_norm}.  As alluded to earlier, we choose the atoms in $\A_{\text{comm}}$ to correspond to distributed controllers implemented on communication graphs that can be constructed by adding a small number of communication links from the set of allowed edges $\mathscr{E}$ to the base QI communication graph $\Gamma_{\text{base}}$.  In order to avoid introducing additional notation we describe the atomic set specified by communication graphs that can be constructed by adding a single communication link from the set $\mathscr{E}$ to the base QI communication graph $\Gamma_{\text{base}}$ -- the presented concepts then extend to the general case in a natural way.  We explain why a controller designer may wish to construct an atomic set specified by more complex communication graphs in \S\ref{sec:discussion}.

\subsubsection*{The atomic set  $\A_{\text{comm}}$}To each communication link $(i,j) \in \mathscr{E}$  we associate the subspace $\mathcal{E}_{ij}$ given by
\begin{equation}
\mathcal{E}_{ij} := \s^\perp(\Gamma_{\text{base}})\cap \s(\Gamma_{\text{base}} + E_{ij}).
\label{eq:link_space}
\end{equation}
\edit{Each subspace $\mathcal{E}_{ij}$ encodes the additional information available to the controller, relative to the base communication graph $\Gamma_{\text{base}}$, that is uniquely due to the added communication link $(i,j)$ from sub-controller $j$ to sub-controller $i$.}  Note that the subspaces $\mathcal{E}_{ij}$ are finite dimensional due to the strong connectedness assumption imposed on $\Gamma_{\text{base}}$, which leads to the equality $\s^\perp(\Gamma_{\text{base}}) = \Y^\perp(\Gamma_{\text{base}})\cap\RHinf^{\leq d(\Gamma_{\text{base}})}$.
\begin{example}
Consider the base QI communication graph $\Gamma_{\text{base}}$ illustrated in Figure \ref{fig:4chain} and specified by \eqref{eq:3-chain_gam}.  This communication graph generates the subspace $\s(\Gamma_{\text{base}})$ shown in \eqref{eq:3-chain_s}.
We consider choosing from two additional links to augment the base communication graph $\Gamma_{\text{base}}$: a directed link from node 1 to node 3, and a directed link from node 3 to node 1.  Then $\mathscr{E} = \{(1,3),(3,1)\}$ and the corresponding subspaces $\mathcal{E}_{ij}$ are given by
\vspace{-4mm}
\begin{equation*}
\begin{array}{rclrcl}
\mathcal{E}_{13} &=& \frac{1}{z^2}\begin{bmatrix} 0 & 0 & 0  \\  0 & 0 & 0 \\ \ast & 0 & 0 \end{bmatrix}, 
\mathcal{E}_{31} &=& \frac{1}{z^2}\begin{bmatrix} 0 & 0 & \ast  \\  0 & 0 & 0 \\ 0 & 0 & 0 \end{bmatrix}.
\end{array}
\vspace{-3mm}
\end{equation*}
\vspace{-3mm}
\hfill{}\QED
\end{example}

The atomic set is then composed of suitably normalized elements of these subspaces:
\begin{equation}
\begin{array}{c}
\edit{\A_{\text{comm}}:= \displaystyle\bigcup_{(i,j) \in \mathscr{E}}\left\{ A \in \mathcal{E}_{ij} \, \big{|} \, \Htwonorm{A} = 1 \right\}.}
\end{array}
\label{eq:comm_set}
\end{equation}
Note that we normalize our atoms relative to the $\Htwo$ norm as this norm is isotropic; hence this normalization ensures that no atom is preferred over another within the family of atoms defined by a subspace $\mathcal{E}_{ij}$.
\edit{The resulting atomic norm, which we denote the \emph{communication link norm}, is defined on elements $X \in \RHinf^{\leq d(\Gamma_\text{base})}$ and  is given by\footnote{\edit{We apply definition \eqref{eq:atom_norm} to the components of $X$ that lie in $\s^\perp(\Gamma_{\text{base}})$ to obtain an atomic norm defined on elements of that space.  We then introduce an unpenalized variable $A_\text{base}\in \Y(\Gamma_{\text{base}})$ to the atomic decomposition so that the resulting penalty function may be applied to elements $X \in \RHinf^{\leq d(\Gamma_\text{base})}$.  The resulting penalty is actually a seminorm on $\RHinf^{\leq d(\Gamma_\text{base})}$ but we refer to it as a norm to maintain consistency with the terminology of \cite{CRPW12}.}}
\begin{equation}
\begin{array}{rl}
\commnorm{X} = \edit{\displaystyle\min_{A_{\text{base}},\left\{A_{ij}\right\}\in\RHinf^{\leq d(\Gamma_\text{base})}}} & \displaystyle\sum_{(i,j)\in\mathscr{E}}\Htwonorm{A_{ij}} \\
\text{s.t.} & X = A_{\text{base}} +\displaystyle \sum_{(i,j) \in \mathscr{E}}A_{ij} \\
&A_{\text{base}} \in \Y(\Gamma_{\text{base}}) \\
& A_{ij} \in \mathcal{E}_{ij} \ \forall (i,j) \in \mathscr{E},
\end{array}
\label{eq:comm_norm}
\end{equation}
when this optimization problem is feasible --  when it is not, we set $\commnorm{X} = \infty$.}
Applying definition \eqref{eq:comm_norm} of the communication link norm to the regularized optimization problem \eqref{eq:reg_opt} yields the convex optimization problem
\begin{equation}
\begin{array}{rl}
 \edit{\minimize{V,A_{\text{base}},\left\{A_{ij}\right\}\in\RHinf^{\leq d(\Gamma_\text{base})}}} & \Htwonorm{\mathfrak{L}(V)}^2 + \lambda \left(\displaystyle\sum_{(i,j)\in\mathscr{E}}\Htwonorm{A_{ij}} \right)\\
\text{s.t.} & \mathfrak{C}(V) = A_{\text{base}} + \displaystyle \sum_{(i,j) \in \mathscr{E}} A_{ij}\\
& A_{\text{base}} \in \Y(\Gamma_{\text{base}}) \\
& A_{ij} \in \mathcal{E}_{ij} \ \forall (i,j) \in \mathscr{E}.
\end{array}
\label{eq:long_opt}
\end{equation}

\edit{Recall that in optimization problem \eqref{eq:reg_opt} our approach to communication architecture design is to induce structure in the term $\mathfrak{C}(V)$ through the use of the communication link norm as a penalty function.  Letting  $\left(V^\star,\{A^\star_{ij}\},A^\star_{\text{base}}\right)$ denote the solution to the optimization problem \eqref{eq:long_opt}, we have that each nonzero $A^\star_{ij}$ in the atomic decomposition of $\mathfrak{C}(V)$ corresponds to an additional link from sub-controller $j$ to sub-controller $i$ being added to the base QI communication graph (in what follows we make precise how the structure of $\mathfrak{C}(V^\star)$ can be used to specify a communication graph). As desired, the communication link norm \eqref{eq:comm_norm} penalizes the use of such additional links, and optimization problem \eqref{eq:long_opt} allows for a tradeoff between communication complexity (as measured by $\sum_{(i,j)\in\mathscr{E}}\Htwonorm{A_{ij}}$) and closed loop performance (as measured by $\Htwonorm{\mathfrak{L}(V)}^2$) of the resulting distributed controller through the regularization weight $\lambda$. Note further that $A^\star_{\text{base}}$ is not penalized by the communication link norm, ensuring that the communication graph defined by the structure of $\mathfrak{C}(V^\star)$ has $\Gamma_{\text{base}}$ as a subgraph. }

\begin{rem}
Optimization problem \eqref{eq:long_opt} is finite dimensional, and hence can be formulated as a second order cone program by associating the finite impulse response transfer matrices $(V,A_{\text{base}},\left\{A_{ij}\right\})$, $ \mathfrak{C}(V)$ and $\mathfrak{L}(V)$ with their matrix representations.  To see this, note that $\Y(\Gamma_{\text{base}})\subseteq\RHinf^{\leq d(\Gamma_{\text{base}})}$, and that by the discussion after the definition \eqref{eq:link_space} of the subspaces $\mathcal{E}_{ij}$, they too satisfy $\mathcal{E}_{ij}\subseteq\RHinf^{\leq d(\Gamma_{\text{base}})}$. Thus the horizon  $d(\Gamma_{\text{base}})$ over which the optimization problem \eqref{eq:long_opt} is solved is finite.
\label{rem:fin_dim}
\end{rem}
\subsection{Co-Design Algorithm and Solution Properties}
\label{sec:algo}
In this section we formally define the communication architecture/control law co-design algorithm in terms of the optimization problem \eqref{eq:long_opt}, and show that it can be used to co-design a strongly connected communication graph $\Gamma$ that generates a QI subspace $\s(\Gamma)$ as defined in \eqref{eq:s_gam}. 

\begin{algorithm}
 \Input{regularization weight $\lambda$, generalized plant $G$, base QI communication graph $\Gamma_{\text{base}}$, edge set $\mathscr{E}$;}
 \Output{designed communication graph adjacency matrix $\Gamma_{\text{des}}$, optimal Youla parameter $Q^\star_{\text{des}}\in \s(\Gamma_{\text{des}})$; }
 \KwInitialize{$\Gamma_{\text{des}} \leftarrow \Gamma_{\text{base}}$, $Q^\star_{\text{des}} \leftarrow 0$;}
 
 \KwGenGraph{
 $\left(V^\star,\{A^\star_{ij}\},A^\star_{\text{base}}\right) \leftarrow$ solution to optimization problem \eqref{eq:long_opt} with regularization weight $\lambda$;
 
 \ForEach(){$(i,j)\in\mathscr{E}$ s.t. $A^\star_{ij}\neq 0$}{$\Gamma_{\text{des}} \leftarrow \Gamma_{\text{des}} + E_{ij}; $}
}
\KwGetK{
$Q^\star_{\text{des}} \leftarrow$ solution to optimization problem \eqref{eq:quad_prog} with distributed constraint $\Y(\Gamma_{\text{des}})$, as specified by Theorem \ref{thm:quad};
}
\Return{$\Gamma_{\text{des}}$, $Q^\star_{\text{des}}$;}
 \caption{Communication Architecture Co-Design}
 \label{lag:co-design}
\end{algorithm}

The co-design procedure is described in Algorithm \ref{lag:co-design}.  The algorithm consists of first solving the regularized optimization problem \eqref{eq:long_opt} to obtain solutions $\left(V^\star,\{A^\star_{ij}\},A^\star_{\text{base}}\right)$.  Using these solutions, we produce the designed communication graph $\Gamma_{\text{des}}$ by augmenting the base QI communication graph $\Gamma_{\text{base}}$ with all edges $(i,j)$ such that $A^\star_{ij} \neq 0$.  In particular, each non-zero term $A^\star_{ij}$ corresponds to an additional edge $(i,j) \in \mathscr{E}$ that the co-designed distributed control law will use -- thus by varying the regularization weight $\lambda$ the controller designer can control how much the use of an additional link is penalized by the optimization problem \eqref{eq:long_opt}.  As $\supp{\Gamma_{\text{base}}}\subseteq\supp{\Gamma_{\text{des}}}\subseteq\supp{\Gamma_{\text{max}}}$ by construction, the designed communication graph $\Gamma_{\text{des}}$ satisfies the assumptions of Proposition \ref{prop:space} -- it is therefore strongly connected, can be physically built, and generates a subspace $\s(\Gamma_{\text{des}})$, according to \eqref{eq:s_gam}, that is QI with respect to $G_{22}$ and that admits a decomposition of the form \eqref{eq:s_decomp}.  
The subspace $\s(\Gamma_{\text{des}})$ thus satisfies the assumptions of Theorem \ref{thm:quad}, meaning that the distributed optimal controller $K^\star_{\text{des}}$ restricted to the designed subspace $\s({\Gamma_{\text{des}}})$ is specified in terms of the solution $Q^\star_{\text{des}}$ to the convex quadratic program \eqref{eq:quad_prog}.  In this way the optimal distributed controller restricted to the designed communication architecture, as well as the performance that it achieves, can be computed exactly.

\edit{Although the solution $V^\star$ to optimization problem \eqref{eq:long_opt} could be used to generate a distributed controller that can be implemented on the designed communication graph $\Gamma_{\text{des}}$, we claim that it is preferable to use the solution $Q^\star_{\text{des}}$ to the non-regularized optimization problem \eqref{eq:quad_prog}.  First, the use of the communication link norm penalty in the optimization problem  \eqref{eq:quad_prog} has the effect of shrinking the solution towards the origin.  This means that the resulting controller specified by $V^\star$ is less aggressive, i.e., has smaller control gains, than the controller specified by the solution to the optimization problem \eqref{eq:quad_prog} with subspace constraint $\Y(\Gamma_{\text{des}})$.}

Second, notice that for two graphs $\Gamma_{ij}$ and $\Gamma_{kl}$ obtained by augmenting the base QI communication graph $\Gamma_{\text{base}}$ with the communication links $(i,j)$ and $(k,l)$, respectively, it holds that $\s(\Gamma_{ij}) + \s(\Gamma_{kl}) \subseteq \s(\supp{\Gamma_{ij} +\Gamma_{kl}})$, with the inclusion being strict in general.
In words, the linear superposition of the subspaces \eqref{eq:s_gam} generated by the two communications graphs $\Gamma_{ij}$ and $\Gamma_{kl}$ is in general a strict subset of the subspace generated by the single communication graph $\supp{\Gamma_{ij}+\Gamma_{kl}}$.  Suppose now that the corresponding solutions $A_{ij}^\star$ and $A^\star_{kl}$ to optimization problem \eqref{eq:long_opt} are non-zero: then $\Gamma_{\text{des}}=\Gamma_{\text{base}} + E_{ij} + E_{kl}$, but the expression $\mathfrak{C}(V^\star)$ lies in the subspace given by $\s(\Gamma_{ij}) + \s(\Gamma_{kl})$.
By the previous discussion $\s(\Gamma_{ij}) + \s(\Gamma_{kl})\subset\s(\Gamma_{\text{des}})$, and thus we are imposing additional structure on the the expression $\mathfrak{C}(V^\star)$ relative to that imposed on the solution to the non-regularized optimization problem \eqref{eq:quad_prog} with subspace constraint $\Y(\Gamma_{\text{des}})$.  This can be interpreted as the controller specified by the structure of $\mathfrak{C}(V^\star)$ not utilizing paths in the communication graph that contain both links $(i,j)$ and $(k,l)$.  These sources of conservatism in the control law are however completely removed if one uses the solution $Q^\star_{\text{des}}$ to the non-regularized optimization problem \eqref{eq:quad_prog}. 

 Thus we have met our objective of developing a convex optimization based procedure for co-designing a distributed optimal controller and the communication architecture upon which it is implemented.  In the next section we discuss the computational complexity of the proposed method and illustrate its efficacy on numerical examples.

\section{Computational Examples}
\label{sec:examples}
We show that the number of scalar optimization variables needed to formulate the regularized optimization problem \eqref{eq:long_opt} scales, up to constant factors, in a manner identical to the number of variables needed to formulate the non-regularized optimization problem \eqref{eq:quad_prog}.  We then illustrate the usefulness of our approach via two examples.

\subsubsection*{Computational Complexity}
We assume that the number of control inputs $p_2$ and the number of measurements $q_2$ scale as $O(n)$, where $n$ is the number of sub-controllers in the system, i.e., we assume that there is an order constant number of actuators and sensors at each sub-controller.  For an element $V \in \RHinf^{\leq d}$,  each term $V^{(t)}$ in its power-series expansion is a real matrix of dimension $O(n) \times O(n)$, and thus $V$ is defined by $O(n^2d)$ scalar variables.  The convex quadratic program \eqref{eq:quad_prog} is therefore specified in terms of $O(n^2d)$ variables.

To describe the number of scalar optimization variables in the regularized optimization problem \eqref{eq:long_opt}, we need to take into account the contributions from $V$, $A_{\text{base}}$ and $\{A_{ij}\}$.  As per the discussion in the previous paragraph, $V$ and $A_{\text{base}}$ are composed of at most $O(n^2d)$ scalar optimization variables.  It can be checked that each $A_{ij}$ has $O(d)$ optimization variables, and hence the collection $\{A_{ij}\}$ contributes $O(d|\mathscr{E}|)$ scalar optimization variables.  Each sub-controller can have at most $O(n)$ additional links originating from it, and thus $|\mathscr{E}|$ scales, at worst, as $O(n^2)$.  It follows that the regularized optimization problem \eqref{eq:long_opt} can also be specified in terms of $O(n^2d)$ scalar optimization variables.  

Finally, we note that the regularized optimization problem \eqref{eq:long_opt} is a second order cone program (SOCP) with at most $O(n^2d)$ second order constraints.  It therefore enjoys favorable iteration complexity that scales as $O(\sqrt{d}n)$ \cite{SOCP_AG}, and its per-iteration complexity is at worst $O(d^3n^6)$ \cite{SOCP_Boyd}, but is typically much less when structure is exploited.  In particular it is not atypical to solve a SOCP with tens to hundreds of thousands of variables \cite{SOCP_benchmark}: noting that $d$ scales at worst as $O(n)$, we therefore expect our method to be applicable to problems with hundreds of sub-controllers.  Further, as we illustrate in the 20 sub-controller ring example below, the computational benefits of our approach compared to a brute force search are already tangible for systems with tens of sub-controllers.

\subsubsection*{6 sub-controller chain system}
Consider a generalized plant \eqref{eq:gen_plant} specified by a tridiagonal matrix $A_{\text{6-chain}}\in\R^{6\times6}$ with randomly generated nonzero entries, $B_2 = C_2 = I_6$, $B_1=C_1^\top = \begin{bmatrix} I_6 & 0_6 \end{bmatrix}$ and $D_{21} = D_{12}^\top = \begin{bmatrix} 0_6 & I_6 \end{bmatrix}$.  The physical topology of the plant $G_{22}$ is that of a 6 subsystem chain (a 3 subsystem chain is illustrated in Figure \ref{fig:4chain}), and therefore the base QI communication graph $\Gamma_{\text{6-chain}}= \bsupp{A_{\text{6-chain}}}$ also defines a 6 sub-controller chain.  We define the set of edges that can be added to the base graph to be
\begin{equation}
\mathscr{E} = \{ (i,j) \in \N\times\N \, \big{|} \, |i-j| = 2\}, \label{eq:example_E} \end{equation} 
i.e., the communication graph/control law co-design task consists of determining which additional directed communication links between second neighbors should be added to the base QI communication graph $\Gamma_{\text{6-chain}}$ to best improve the performance of the distributed optimal controller implemented on the resulting augmented communication graph.  

\begin{figure}[h]
\centering
\includegraphics[width=2.7in]{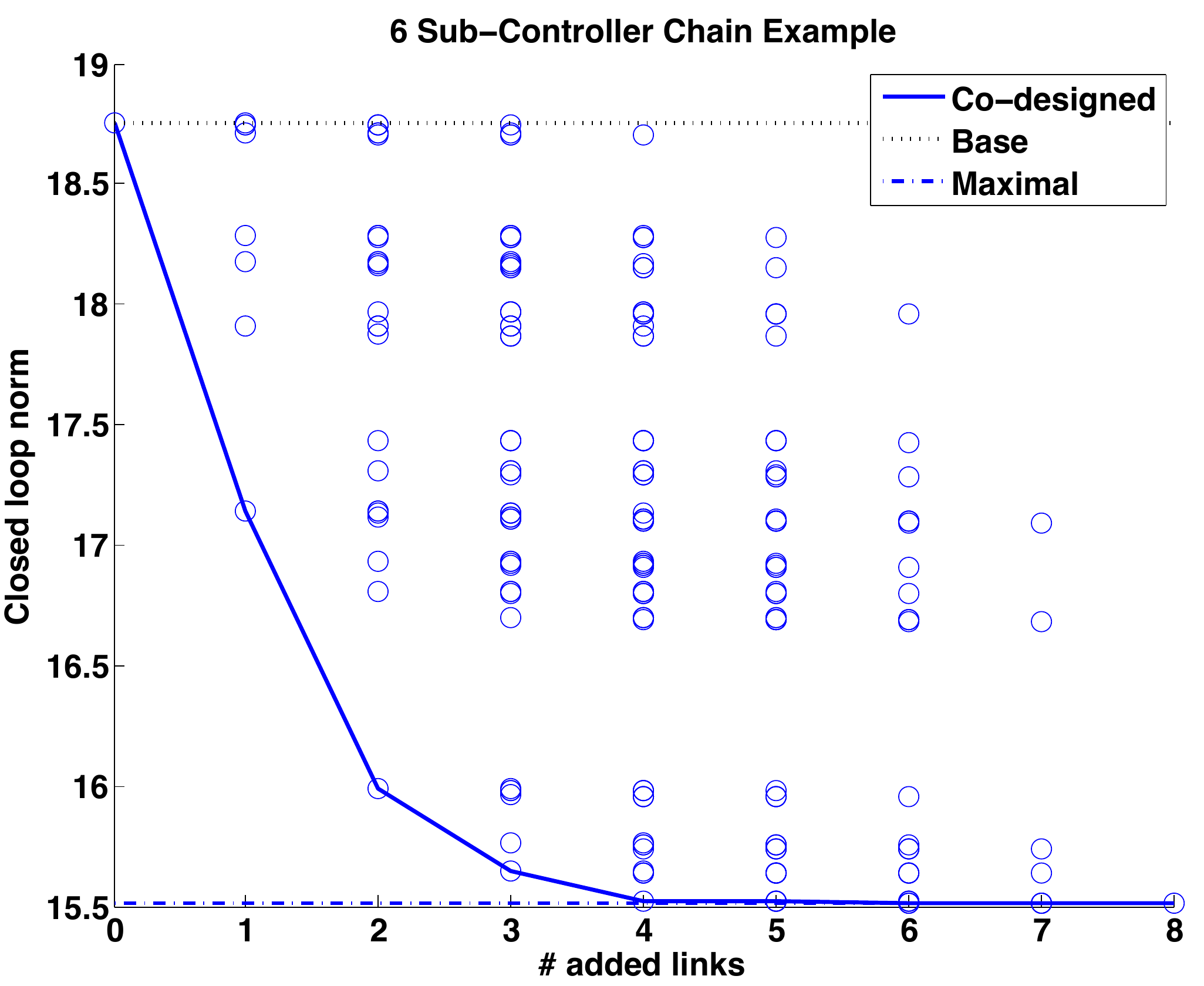}
\caption{\small The closed loop norms achieved by distributed optimal controllers implemented on communication graphs constructed by adding $k=1,\dots,|\mathscr{E}|$ links to the base QI communication graph  $\Gamma_{\text{6-chain}}$ are plotted as blue circles.  The solid blue line denotes the performance achieved by distributed optimal controllers implemented on the communication graphs identified by the co-design procedure described in Algorithm \ref{lag:co-design}.  The dotted and dashed lines indicate the closed loop norm achieved by the distributed optimal controllers implemented on the base and maximal QI communication graphs, respectively.}
\label{fig:6chaingraph_space}
\end{figure}

In order to assess the efficacy of the proposed method in uncovering communication topologies that are well suited to distributed optimal control, we first computed the optimal closed loop performance achievable by a distributed controller implemented on every possible communication graph that can be constructed by augmenting the base QI communicating graph $\Gamma_{\text{6-chain}}$ with $k=1,\dots,|\mathscr{E}|$ additional links drawn from the set $\mathscr{E}$.   In particular, we exhaustively explored the QI communication graph set $\mathscr{G}$ and computed the achievable closed loop norms -- these closed loop norms are plotted as blue circles in Figure \ref{fig:6chaingraph_space}.
 We then performed the co-design procedure described in Algorithm \ref{lag:co-design} for different values of regularization weight $\lambda \in [0,50]$.  The resulting closed loop norms achieved by the co-designed communication architecture/control law are plotted as a solid blue line in Figure \ref{fig:6chaingraph_space}.  We also plot the closed loop norms achieved by controllers implemented using the base and maximal QI communication graphs.
 
We observe that as the regularization weight $\lambda$ is increased, simpler communication topologies are generated by the co-design procedure.  Further, our algorithm is able to successfully identify the optimal communication topology and the corresponding distributed optimal control law for every fixed number of additional links.

\subsubsection*{20 sub-controller ring system}
\begin{figure}[h]
\centering
\includegraphics[width=2.7in]{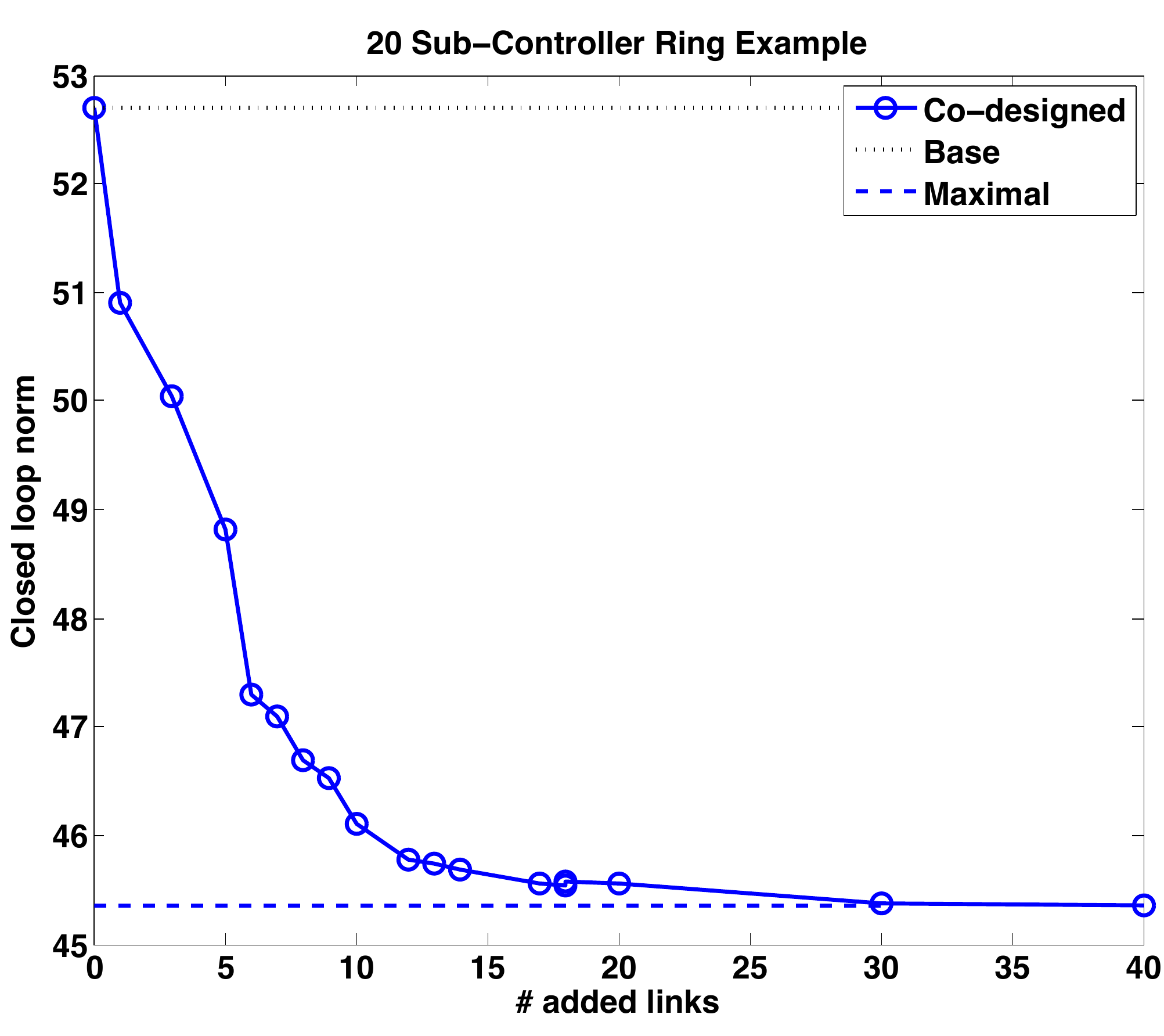}
\caption{ \small The solid blue line denotes the performance achieved by distributed optimal controllers implemented on the communication graphs identified by the co-design procedure described in Algorithm \ref{lag:co-design}.  The dotted and dashed lines indicate the closed loop norm achieved by the distributed optimal controllers implemented on the base and maximal QI communication graphs, respectively.}
\label{fig:20chain}
\end{figure}

Consider a generalized plant \eqref{eq:gen_plant} specified by a matrix $A_{\text{20-ring}}\in\R^{20\times20}$ with $(i,j)$th entry set to a nonzero randomly generated number if $|i-j|\leq 1$ where the subtraction is modulo 20 \edit{(e.g., 1-20 = 1)}, and 0 otherwise.  The additional state-space parameters are given by $B_2 = C_2 = I_{20}$, $B_1=C_1^\top = \begin{bmatrix} I_{20} & 0_{20} \end{bmatrix}$ and $D_{21} = D_{12}^\top = \begin{bmatrix} 0_{20}& I_{20} \end{bmatrix}$.  For the example considered below, $|\lambda_{\max}(A_{\text{20-ring}})|=2.91.$
The physical topology of the plant $G_{22}$ is that of a 20 subsystem ring, i.e., a chain topology with first and last nodes connected, and therefore the base QI communication graph $\Gamma_{\text{20-ring}}= \bsupp{A_{\text{20-ring}}}$ also defines a 20 sub-controller ring.  We again define the set of edges $\mathscr{E}$ that can be added to the base graph to be those between second neighbors as in \eqref{eq:example_E}.  In this case, the QI communication graph set $
\mathscr{G}$ is too large to exhaustively explore: in particular \edit{$|\mathscr{G}| = 2^{40}\approx 10^{12}$}.  We performed the co-design procedure described in Algorithm \ref{lag:co-design} for different values of regularization weight $\lambda \in [0,1000]$.  The resulting closed loop norms achieved by the co-designed communication architecture/control law are plotted as a solid blue line in Figure \ref{fig:20chain}.  We also plot the closed loop norms achieved by controllers implemented using the base and maximal QI communication graphs.
 We observe again that as the regularization weight $\lambda$ is increased, simpler and simpler communication topologies are designed.  Notice that our method selected $10$ carefully placed communication links to add to the base QI communication graph, leading to a closed loop performance only 2\% higher than that achieved by the optimal controller implemented using the maximal QI communication graph.

\section{Discussion}
\label{sec:discussion}
\noindent\textbf{Optimal structural recovery:} It is shown in \cite{MC_RFD_TAC} that the variational solution to an $\Htwo$ optimal control problem augmented with an atomic norm that penalizes the use of actuators can succeed in identifying an optimal actuation architecture when the dynamics of the plant satisfy certain conditions.  The numerical experiments of \S\ref{sec:examples} provide empirical evidence that our approach to communication architecture design identifies optimally structured controllers as well -- it is of interest to see whether conditions analogous to those of \cite{MC_RFD_TAC} can provide theoretical support to the empirical success of our approach.  

\noindent\textbf{The $k$-communication link norm:} The communication link norm was defined in terms of atoms corresponding to communication graphs constructed by adding a single link to the base QI communication graph.  However it is possible to include atoms corresponding to communication graphs augmented with at most $k$-links instead, for any positive integer $k$; denote the resulting $k$-communication link norm by $\kcommnorm{k}{\cdot}$.  If the atoms are suitably normalized,\footnote{In particular, elements $A \in \mathscr{A}_{\text{k-comm}}$ constrained to lie in a subspace $\mathcal{E}$ should be normalized as $\Htwonorm{A} = \left(\text{card}\left(\mathcal{E}\right) + \kappa\right)^{-\frac{1}{2}}$, where $\kappa>0$ is a positive constant that controls how much a single atom of larger cardinality is preferred over several atoms of lower cardinality.} for all positive integers  $k_1$ and $k_2$ satisfying $k_1\leq k_2$
it then holds that $\kcommnorm{k_1}{G} \leq \kcommnorm{k_2}{G}$ for all transfer matrices $G$ satisfying $\kcommnorm{k_1}{G}<\infty$.  Geometrically, restricted to the domain of $\kcommnorm{k_1}{\cdot}$, the unit ball of $\kcommnorm{k_2}{\cdot}$ is an inner approximation to that of $\kcommnorm{k_1}{\cdot}$, and may therefore lead to simpler communication graphs when used as a penalty function in the regularized optimization problem \eqref{eq:reg_opt}.  How to choose $k$ will presumably be informed by the aforementioned conditions on optimal communication structure recovery, as well as by computational considerations, as the number of elements $\{A_{ij}\}$ required to implement the $k$-communciation link norm scales as $O(n^{2k})$.

\noindent\textbf{Constructing base QI communication graphs:}  \edit{The structural assumptions made on $(A,B_2,C_2)$ in \S\ref{sec:delays} are needed to ensure that the base QI communication graph as specified in \eqref{eq:gam_min} is strongly connected and generates a QI subspace.  However, as we note in Remark \ref{rem:base}, any strongly connected communication topology leading to a QI subspace can be used as the base QI communication graph.}
 Exploring how to construct base QI communication graphs in a principled way when the structural assumptions on $(A,B_2,C_2)$ are relaxed, perhaps utilizing the methods in \cite{RM12}, is an interesting direction for future work.  We emphasize however that the rest of the discussion in \S\ref{sec:delays} remains valid once a base QI communication graph is identified even if the structural assumptions on $(A,B_2,C_2)$ are relaxed .  We also note that these issues are a consequence of the communication protocol imposed between sub-controllers -- determining alternative communication protocols that allow the structural assumptions to be relaxed is also an interesting direction for future work.  

\noindent\textbf{Scalability:} Although we expect the methods presented to be applicable to systems composed of hundreds of sub-controllers, it is important that the general approach of the RFD framework be applicable to truly large-scale systems composed of heterogeneous subsystems.  The limits on the scalability of our proposed method are due to the underlying controller synthesis method \cite{LD14}, as opposed to being inherent to the communication link norm.  To that end we have been pursuing \emph{localized optimal control} \cite{WMD_CDC14} as a scalable distributed optimal controller synthesis method -- an interesting direction for future work will be to see if communication architecture co-design can be incorporated into the localized optimal control framework.
%

\subsubsection*{Acknowledgements}
The author would like to thank V. Chandrasekaran, A. Lamperski and J. C. Doyle for insightful discussions, V. D. Jonsson for her careful revisions, K. Dvijotham for references on the complexity of SOCPs, and the anonymous reviewers for their helpful comments.
\bibliographystyle{IEEEtran}
\bibliography{../../biblio/decentralized,../../biblio/comms,../../biblio/regularization,../../biblio/matni,../../biblio/sys_id}

\begin{thebibliography}{10}
\providecommand{\url}[1]{#1}
\csname url@rmstyle\endcsname
\providecommand{\newblock}{\relax}
\providecommand{\bibinfo}[2]{#2}
\providecommand\BIBentrySTDinterwordspacing{\spaceskip=0pt\relax}
\providecommand\BIBentryALTinterwordstretchfactor{4}
\providecommand\BIBentryALTinterwordspacing{\spaceskip=\fontdimen2\font plus
\BIBentryALTinterwordstretchfactor\fontdimen3\font minus
  \fontdimen4\font\relax}
\providecommand\BIBforeignlanguage[2]{{%
\expandafter\ifx\csname l@#1\endcsname\relax
\typeout{** WARNING: IEEEtran.bst: No hyphenation pattern has been}%
\typeout{** loaded for the language `#1'. Using the pattern for}%
\typeout{** the default language instead.}%
\else
\language=\csname l@#1\endcsname
\fi
#2}}

\bibitem{M_CDC13_codesign}
N.~Matni, ``Communication delay co-design in $\mathcal{H}_2$ decentralized
  control using atomic norm minimization,'' in \emph{Decision and Control
  (CDC), 2013 IEEE 52nd Annual Conference on}, Dec 2013, pp. 6522--6529.

\bibitem{MC_RFD_TAC}
\BIBentryALTinterwordspacing
N.~Matni and V.~Chandrasekaran, ``Regularization for design,'' \emph{Submitted
  to the IEEE Transactions on Automatic Control}, 2015. [Online]. Available:
  \url{http://arxiv.org/abs/1404.1972}
\BIBentrySTDinterwordspacing

\bibitem{CRPW12}
V.~Chandrasekaran, B.~Recht, P.~Parrilo, and A.~Willsky,
  ``\BIBforeignlanguage{English}{The convex geometry of linear inverse
  problems},'' \emph{\BIBforeignlanguage{English}{Foundations of Computational
  Mathematics}}, vol.~12, pp. 805--849, 2012.

\bibitem{Bon91}
F.~Bonsall, ``A general atomic decomposition theorem and {Banach}'s closed
  range theorem,'' \emph{The Quarterly Journal of Mathematics}, vol.~42, no.~1,
  pp. 9--14, 1991.

\bibitem{Pis86}
G.~Pisier, \emph{Probabilistic methods in the geometry of {Banach}
  spaces}.\hskip 1em plus 0.5em minus 0.4em\relax Springer, 1986.

\bibitem{W68}
H.~S. Witsenhausen, ``A counterexample in stochastic optimum control,''
  \emph{SIAM Journal on Control}, vol.~6, no.~1, pp. 131--147, 1968.

\bibitem{RL06}
M.~Rotkowitz and S.~Lall, ``A characterization of convex problems in
  decentralized control,'' \emph{Automatic Control, IEEE Transactions on},
  vol.~51, no.~2, pp. 274--286, 2006.

\bibitem{RCL10}
M.~Rotkowitz, R.~Cogill, and S.~Lall, ``Convexity of optimal control over
  networks with delays and arbitrary topology,'' \emph{International Journal of
  Systems, Control and Communications}, vol.~2, no. 1-3, pp. 30--54, 2010.

\bibitem{LL11_QI}
L.~Lessard and S.~Lall, ``Quadratic invariance is necessary and sufficient for
  convexity,'' in \emph{American Control Conference (ACC), 2011}, June 2011,
  pp. 5360--5362.

\bibitem{BV05}
B.~Bamieh and P.~G. Voulgaris, ``A convex characterization of distributed
  control problems in spatially invariant systems with communication
  constraints,'' \emph{Sys. \& Control Letters}, vol.~54, no.~6, pp. 575 --
  583, 2005.

\bibitem{MMRY12}
A.~Mahajan, N.~Martins, M.~Rotkowitz, and S.~Yuksel, ``Information structures
  in optimal decentralized control,'' in \emph{Decision and Control (CDC), 2012
  IEEE 51st Annual Conference on}, Dec 2012, pp. 1291--1306.

\bibitem{LDXX}
A.~Lamperski and J.~Doyle, ``Output feedback $\mathcal{H}_2$ model matching for
  decentralized systems with delays,'' in \emph{American Control Conference
  (ACC), 2013}, June 2013, pp. 5778--5783.

\bibitem{LD14}
A.~Lamperski and J.~C. Doyle, ``The $\mathcal{H}_2$ control problem for
  quadratically invariant systems with delays,'' \emph{Automatic Control, IEEE
  Transactions on. To appear.}, 2015.

\bibitem{SBTR12}
P.~Shah, B.~N. Bhaskar, G.~Tang, and B.~Recht, ``Linear system identification
  via atomic norm regularization,'' in \emph{Decision and Control (CDC), 2012
  IEEE 51st Annual Conference on}, 2012, pp. 6265--6270.

\bibitem{Fazel}
M.~Fazel, H.~Hindi, and S.~Boyd, ``A rank minimization heuristic with
  application to minimum order system approximation,'' in \emph{American
  Control Conference (ACC), 2001}, June 2001, pp. 4734--4739.

\bibitem{MR_CDC14}
\BIBentryALTinterwordspacing
N.~Matni and A.~Rantzer, ``Low-rank and low-order decompositions for local
  system identification,'' \emph{CoRR}, vol. arXiv:1403.7175, 2014. [Online].
  Available: \url{http://arxiv.org/abs/1403.7175}
\BIBentrySTDinterwordspacing

\bibitem{LjungNewOld}
L.~Ljung, ``Some classical and some new ideas for identification of linear
  systems,'' \emph{Journal of Control, Automation and Electrical Systems},
  vol.~24, no. 1-2, pp. 3--10, 2013.

\bibitem{FLJ11}
M.~Fardad, F.~Lin, and M.~Jovanovic, ``Sparsity-promoting optimal control for a
  class of distributed systems,'' in \emph{American Control Conference (ACC),
  2011}, June 2011, pp. 2050--2055.

\bibitem{LFJ13}
F.~Lin, M.~Fardad, and M.~R. Jovanovic, ``Design of optimal sparse feedback
  gains via the alternating direction method of multipliers,'' \emph{Automatic
  Control, IEEE Transactions on}, vol.~58, no.~9, pp. 2426--2431, 2013.

\bibitem{JRM14}
V.~Jonsson, A.~Rantzer, and R.~Murray, ``A scalable formulation for engineering
  combination therapies for evolutionary dynamics of disease,'' in
  \emph{American Control Conference (ACC), 2014}, June 2014, pp. 2771--2778.

\bibitem{DLFM12}
N.~Dhingra, F.~Lin, M.~Fardad, and M.~R. Jovanovic, ``On identifying sparse
  representations of consensus networks,'' in \emph{3rd IFAC Workshop on
  Distributed Estimation and Control in Networked Systems, Santa Barbara, CA},
  2012, pp. 305--310.

\bibitem{XB07}
L.~Xiao, S.~Boyd, and S.-J. Kim, ``Distributed average consensus with
  least-mean-square deviation,'' \emph{J. Parallel Distrib. Comput.}, vol.~67,
  no.~1, pp. 33--46, Jan. 2007.

\bibitem{FLJ13}
M.~Fardad, F.~Lin, and M.~Jovanovic, ``Design of optimal sparse interconnection
  graphs for synchronization of oscillator networks,'' \emph{Automatic Control,
  IEEE Transactions on}, vol.~59, no.~9, pp. 2457--2462, Sept 2014.

\bibitem{MC_CDC14}
N.~Matni and V.~Chandrasekaran, ``Regularization for design,'' in
  \emph{Decision and Control (CDC), 2014 IEEE 53rd Annual Conference on}, Dec
  2014, pp. 1111--1118.

\bibitem{DJL14}
N.~Dhingra, M.~Jovanovic, and Z.-Q. Luo, ``An {ADMM} algorithm for optimal
  sensor and actuator selection,'' in \emph{Decision and Control (CDC), 2014
  IEEE 53rd Annual Conference on}, Dec 2014, pp. 4039--4044.

\bibitem{Pol13}
B.~Polyak, M.~Khlebnikov, and P.~Shcherbakov, ``An {LMI} approach to structured
  sparse feedback design in linear control systems,'' in \emph{Control
  Conference (ECC), 2013 European}, July 2013, pp. 833--838.

\bibitem{ZDG96}
K.~Zhou, J.~C. Doyle, and K.~Glover, \emph{Robust and Optimal Control}.\hskip
  1em plus 0.5em minus 0.4em\relax Upper Saddle River, NJ, USA: Prentice-Hall,
  Inc., 1996.

\bibitem{GR01_AGT}
C.~D. Godsil, G.~Royle, and C.~Godsil, \emph{Algebraic Graph Theory}.\hskip 1em
  plus 0.5em minus 0.4em\relax Springer New York, 2001, vol. 207.

\bibitem{CDS98}
S.~S. Chen, D.~L. Donoho, and M.~A. Saunders, ``Atomic decomposition by basis
  pursuit,'' \emph{SIAM Rev.}, vol.~43, no.~1, pp. 129--159, Jan. 2001.

\bibitem{CRT06}
E.~Candes, J.~Romberg, and T.~Tao, ``Robust uncertainty principles: exact
  signal reconstruction from highly incomplete frequency information,''
  \emph{Info. Theory, IEEE Trans. on}, vol.~52, no.~2, pp. 489--509, Feb. 2006.

\bibitem{Donoho04}
D.~L. Donoho, ``For most large underdetermined systems of linear equations the
  minimal $\ell_1$-norm solution is also the sparsest solution,'' \emph{Comm.
  Pure Appl. Math}, vol.~59, pp. 797--829, 2004.

\bibitem{fazelThesis}
M.~Fazel, ``Matrix rank minimization with applications,'' Ph.D. dissertation,
  PhD thesis, Stanford University, 2002.

\bibitem{RFP10}
B.~Recht, M.~Fazel, and P.~A. Parrilo, ``Guaranteed minimum-rank solutions of
  linear matrix equations via nuclear norm minimization,'' \emph{SIAM Rev.},
  vol.~52, no.~3, pp. 471--501, Aug. 2010.

\bibitem{CR12}
E.~J. Candes and B.~Recht, ``Exact matrix completion via convex optimization,''
  \emph{Found. Comp. Math.}, vol.~9, no.~6, pp. 717--772, Dec. 2009.

\bibitem{SOCP_AG}
F.~Alizadeh and D.~Goldfarb, ``Second-order cone programming,''
  \emph{Mathematical programming}, vol.~95, no.~1, pp. 3--51, 2003.

\bibitem{SOCP_Boyd}
M.~S. Lobo, L.~Vandenberghe, S.~Boyd, and H.~Lebret, ``Applications of
  second-order cone programming,'' \emph{Linear Algebra and its Applications},
  vol. 284, no. 1--3, pp. 193 -- 228, 1998.

\bibitem{SOCP_benchmark}
\BIBentryALTinterwordspacing
H.~D. Mittelmann. (2014) {MISOCP} and large {SOCP} benchmark. [Online].
  Available: \url{http://plato.asu.edu/ftp/socp.html}
\BIBentrySTDinterwordspacing

\bibitem{RM12}
M.~Rotkowitz and N.~Martins, ``On the nearest quadratically invariant
  information constraint,'' \emph{Automatic Control, IEEE Transactions on},
  vol.~57, no.~5, pp. 1314--1319, May 2012.

\bibitem{WMD_CDC14}
Y.-S. Wang, N.~Matni, and J.~Doyle, ``Localized {LQR} optimal control,'' in
  \emph{Decision and Control (CDC), 2014 IEEE 53rd Annual Conference on}, Dec
  2014, pp. 1661--1668.

\end{thebibliography}

\end{document}